\documentclass[reqno]{amsart}
\usepackage{setspace,tikz,xcolor,mathrsfs,listings,multicol}
\usepackage{rotating}
\usepackage[vcentermath]{youngtab}
\usepackage{enumerate}
\usepackage[all,cmtip]{xy}
\usetikzlibrary{arrows,matrix}
\usepackage[margin=1.25in]{geometry}

\usepackage[colorlinks=true, pdfstartview=FitV, linkcolor=blue, citecolor=blue, urlcolor=blue]{hyperref}

\newcommand{\g}{\mathfrak{g}}
\newcommand{\fsl}{\mathfrak{sl}} 
\newcommand{\asl}{\widehat{\fsl}} 

\newcommand{\cl}[1]{\overline{#1}}
\newcommand{\clfw}{\cl{\Lambda}} 
\newcommand{\inner}[2]{\left\langle #1, #2 \right\rangle} 
\newcommand{\iso}{\cong}

\newcommand{\exmid}{\hspace{5pt}\vline\hspace{5pt}}

\DeclareMathOperator{\wt}{wt} 
\DeclareMathOperator{\pr}{pr} 

\newcommand{\mcM}{\mathcal{M}}

\newcommand{\ZZ}{\mathbb{Z}}

\newcommand{\CC}{\mathbb{C}}

\definecolor{darkred}{rgb}{0.7,0,0} 
\newcommand{\defn}[1]{{\color{darkred}\emph{#1}}} 

\usepackage{listings}
\lstdefinelanguage{Sage}[]{Python}
{morekeywords={False,sage,True},sensitive=true}
\definecolor{dblackcolor}{rgb}{0.0,0.0,0.0}
\definecolor{dbluecolor}{rgb}{0.01,0.02,0.7}
\definecolor{dgreencolor}{rgb}{0.2,0.4,0.0}
\definecolor{dgraycolor}{rgb}{0.30,0.3,0.30}

\usepackage{xparse}

\makeatletter
\protected\def\specialmergetwolists{%
  \begingroup
  \@ifstar{\def\cnta{1}\@specialmergetwolists}
    {\def\cnta{0}\@specialmergetwolists}%
}
\def\@specialmergetwolists#1#2#3#4{%
  \def\tempa##1##2{%
    \edef##2{%
      \ifnum\cnta=\@ne\else\expandafter\@firstoftwo\fi
      \unexpanded\expandafter{##1}%
    }%
  }%
  \tempa{#2}\tempb\tempa{#3}\tempa
  \def\cnta{0}\def#4{}%
  \foreach \x in \tempb{%
    \xdef\cnta{\the\numexpr\cnta+1}%
    \gdef\cntb{0}%
    \foreach \y in \tempa{%
      \xdef\cntb{\the\numexpr\cntb+1}%
      \ifnum\cntb=\cnta\relax
        \xdef#4{#4\ifx#4\empty\else,\fi\x#1\y}%
        \breakforeach
      \fi
    }%
  }%
  \endgroup
}
\makeatother

\usepackage{xparse}
\DeclareDocumentCommand\rpp{ m m g }{
	\foreach \x [count=\s from 1] in {#1}{
	        {\ifnum\s=1
	                \draw (0,-\s)--(\x,-\s);
	                \fi}
	   \draw (0,-\s-1) to (\x,-\s-1);
	   \foreach \y in {0, ..., \x} {\draw (\y,-\s)--(\y,-\s-1);}
	}
	\specialmergetwolists{/}{#1}{#2}\ziplist
	\foreach \x/\y [count=\yi from 1] in \ziplist{
	    \node[anchor=west,font=\small] at (\x,-\yi - .5) {$\y$};
	}
	\IfValueT {#3}
	{\foreach \z [count=\zi from 1] in {#3} {\node[anchor=east,font=\small] at (0,-\zi - .5) {$\z$};}}
	{}
}

\theoremstyle{plain}
\newtheorem{theorem}{Theorem}[section]
\newtheorem{lemma}[theorem]{Lemma}

\newtheorem{proposition}[theorem]{Proposition}

\theoremstyle{definition}

\newtheorem{example}[theorem]{Example}
\newtheorem{remark}[theorem]{Remark}
\numberwithin{equation}{section}

\usepackage[colorinlistoftodos]{todonotes}

\begin{document}
\title{Kirillov--Reshetikhin crystals $B^{1,s}$ for $\asl_n$ using Nakajima monomials}

\author[E.~Gunawan]{Emily Gunawan}
\address[E.~Gunawan]{
Department of Mathematics, 
University of Connecticut, 
Storrs, CT 06269, USA}
\email{egunawan@math.umn.edu}
\urladdr{http://egunawan.github.io}

\author[T.~Scrimshaw]{Travis Scrimshaw}
\address[T.~Scrimshaw]{School of Mathematics and Physics, The University of Queensland, St. Lucia, QLD 4072, Australia}
\email{tcscrims@gmail.com}
\urladdr{https://sites.google.com/view/tscrim/home}

\keywords{crystal, Nakajima monomial, quantum group, Kirillov--Reshetikhin crystal}
\subjclass[2010]{05E10, 17B37}

\thanks{The authors were partially supported by the National Science Foundation RTG grant NSF/DMS-1148634.}

\begin{abstract}
We give a realization of the Kirillov--Reshetikhin crystal $B^{1,s}$ using Nakajima monomials for $\asl_n$ using the crystal structure given by Kashiwara. We describe the tensor product $\bigotimes_{i=1}^N B^{1,s_i}$ in terms of a shift of indices, allowing us to recover the Kyoto path model. Additionally, we give a model for the KR crystals $B^{r,1}$ using Nakajima monomials.
\end{abstract}

\maketitle

\section{Introduction}
\label{sec:introduction}

A special class of finite-dimensional modules of the derived subalgebra of the (Drinfel'd--Jimbo) quantum group $U_q'(\asl_n)$ called Kirillov--Reshetikhin (KR) modules have received significant attention over the past 20 years. KR modules have many remarkable properties and deep connections with mathematical physics. For example, KR modules arise in the study of certain solvable lattice models~\cite{BBB16,JM95,KP84}. Their characters (resp.\ $q$-characters~\cite{FM01,FR99}) satisfy the Q-system (resp.\ T-system) relations, which come from a certain cluster algebra~\cite{dFK09,Hernandez10,Nakajima03II}. This gives a fermionic formula interpretation and a relation to the string hypothesis in the Bethe ansatz for solving Heisenberg spin chains. In untwisted affine types, the graded characters of (Demazure submodules of) tensor products of certain KR modules, the fundamental representations, are also (nonsymmetric) Macdonald polynomials at $t = 0$~\cite{LNSSS14,LNSSS14II} (\cite{LNSSS15}).

In the seminal papers~\cite{K90,K91}, Kashiwara defined the crystal basis of a representation of a quantum group, which is a basis that is well-behaved in the $q \to 0$ limit and affords a combinatorial description. Furthermore, he showed that every irreducible highest weight representation admits a crystal basis $B(\lambda)$. While KR modules are cyclic modules, they are not highest weight modules. Yet, KR modules for $U_q'(\asl_n)$ admit crystal bases~\cite{KKMMNN92} (conjecturally for all affine types~\cite{HKOTT02,HKOTY99}, which is known for non-exceptional types~\cite{OS08} and some other special cases~\cite{JS10,KMOY07,Yamane98}), which are known as Kirillov--Reshetikhin (KR) crystals, and contain even further connections to mathematical physics. For example, KR crystals are in bijection with combinatorial objects that arise naturally from the Bethe ansatz called rigged configurations~\cite{DS06,KKR86,KR86,KSS02}. KR crystals $B^{1,s}$ can be used to model the Takahashi--Satsuma box-ball system~\cite{TS90}, where rigged configurations are invariants called action-angle variables~\cite{KOSTY06,Takagi05}. They are also perfect crystals~\cite{FOS10}, and therefore, they can be used to construct the Kyoto path model~\cite{KKMMNN91,KKMMNN92,OSS03IV}, which came from the study of 2D solvable lattice models and Baxter's corner transfer matrix~\cite{B89}.

Despite intense study, relatively little is understood about KR crystals. In particular, there is currently not a combinatorial model for KR crystals where all crystal operators are given by the same rules, the model is valid for general $B^{r,s}$, and the model is given uniformly across all affine types. By using the decomposition into $U_q(\fsl_n)$-crystals and the Dynkin diagram automorphism, we can lift the tableau model of~\cite{KN94} to a model for KR crystals for $U_q'(\asl_n)$~\cite{Shimozono02}. However, this process obscures the affine crystal operators as it uses the promotion operator of Sch\"utzenberger~\cite{Sch72}, and so it is desirable to have a model where all of the crystal operators are given by the same rules. A similar procedure was utilized in~\cite{FOS09,JS10}, but using type-dependent information. Yet, it cannot work for type $E_8^{(1)}$ due to the Dynkin diagram not admitting any non-trivial automorphisms.

Partial progress has been made on a type-independent construction of KR crystals. Naito and Sagaki constructed a uniform model across all types for tensor products of the form $\bigotimes_{k=1}^N B^{r_k,1}$ by using the usual crystal structure on Lakshmibai--Seshadri (LS) paths for level-zero representations and projecting onto the classical weight space~\cite{NS03,NS05,NS06II,NS06,NS08II}. There is another description of these paths called quantum LS paths~\cite{LNSSS14,LNSSS16}. Lenart and Lubovsky performed a similar construction using a discrete version of quantum LS paths called the quantum alcove path model~\cite{LL15}. Yet, it is not known how to extend these models for general $B^{r,s}$. On the other side, models for $B^{r,s}$ were constructed in~\cite{Kus13,Kus16,Kwon13} for type $A_n^{(1)}$, but these are not known to extend (uniformly) to other affine types.

There is a $t$-analog of $q$-characters (or $q,t$-characters for short) that was studied by Nakajima~\cite{Nakajima01,Nakajima03II,Nakajima03,Nakajima04,Nakajima10} using quiver varieties to show their existence in simply-laced types. From his construction, Nakajima gave a $U_q(\fsl_n)$-crystal structure on the monomials that appear in a $q$-character~\cite{Nakajima03}. Kashiwara~\cite{K03II} independently constructed a different crystal structure on the $q$-character monomials. These two crystal structures were later simultaneously generalized by Sam and Tingley~\cite{ST14}, where the connection with quiver varieties was expanded. This model for crystals is known as the Nakajima monomial model. We note that Kashiwara's crystal structure can be used to realize the crystal $B(\lambda)$ of a highest weight representation for any symmetrizable Kac--Moody Lie algebra $\g$, but Nakajima's requires the Dynkin diagram of $\g$ to not contain an odd length cycle. Kashiwara's crystal structure was extended to $B(\infty)$, the crystal of the lower half of the quantum group or the Verma module of highest weight $0$, by Kang, Kim, and Shin~\cite{KKS07}.

For an extremal level-zero crystal $B(\lambda)$ (so $\lambda$ is a level-zero weight), there exists an automorphism $\kappa$ such that $B(\lambda) / \kappa \iso \bigotimes_{i=1}^N B^{r_i, 1}$ as $U_q'(\asl_n)$-crystals, where $\lambda = \sum_{i=1}^N \Lambda_{r_i}$. This was the construction of Naito and Sagaki previously mentioned, where the description was given explicitly and uniformly in terms of LS paths. A similar construction was given for Nakajima monomials by Hernandez and Nakajima~\cite{HN06}, where the automorphism was constructed type-by-type and does not include all nodes of type $E_{7,8}^{(1)}$ and $E_6^{(2)}$.


Cluster algebras~\cite{FZ02} also have strong connections to characters of KR crystals and Nakajima monomials. Hernandez and Leclerc gave an algorithm to compute $q$-characters as cluster variables of a cluster algebra from a certain semi-infinite quiver~\cite{HL16}. For a double Bruhat cell $G^{u,v} = BuB \cap B_- v B_-$, the coordinate ring $\CC[G^{u,v}]$ is an upper cluster algebra, and the generalized minors of $G^{u,v}$ are cluster variables in $\CC[G^{u,v}]$~\cite{BFZ05}. Kanakubo and Nakashima showed that the generalized minors of $G^{u,e}$ can be expressed as the sum over the Nakajima monomials in a Demazure subcrystal~\cite{KN15}. A further connection between cluster variables in $\CC[G^{u,u^{-1}}]$, with $u$ a Coxeter element, and representation theory was given by Rupel, Stella, and Williams in~\cite{RSW16}. In particular, they show the regular cluster variables in the coordinate ring of the universal central extension of the loop group of $SL_n$ are restrictions of generalized minors of level-zero representations. The specialization of nonsymmetric Macdonald polynomials at $t = \infty$ can also be described as a graded character of a Demazure submodule of a tensor product of fundamental representations~\cite{NNS16}, which satisfy the quantum Q-system relations of~\cite{dFK15,dFK16}.

The main result of this paper is a model for the KR crystal $B^{1,s}$ in type $A_n^{(1)}$ using Nakajima monomials. An important aspect of this construction is that we do not need to apply a quotient automorphism $\kappa$ or use an extremal level-zero crystal for $(B^{1,1})^{\otimes N}$. From this construction, we are able to describe the tensor product of KR crystals $\bigotimes_{i=1}^N B^{1,s_i}$ using only Nakajima monomials ({\it i.e.}, no tensor products). Furthermore, by using the characterization of $B(\lambda)$ given by Kim~\cite[Cor.~4.9]{Kim05}, we are then able to recover the Kyoto path model. From this construction, we are able to relate the models of~\cite{ST14,Tingley08} with the Kyoto path model through the Nakajima monomial model.
Additionally, we give some evidence that our construction can be viewed as a crystal-theoretic interpretation of the purely algebraic reformulation of $q,t$-characters by Hernandez~\cite{Hernandez04} using a $t$-analog of screening operators, as well as the fusion construction of~\cite{KKMMNN91,KKMMNN92}.

Our results can be considered as further evidence of a deep connection between cluster algebras, $q,t$-characters, and KR crystals. Indeed, KR modules can be considered as representations of the loop algebra $\fsl_n[t,t^{-1}]$ and are closely related to level-zero representations as mentioned above. Moreover, after removing certain $0$-arrows, KR crystals are Demazure subcrystals of affine highest weight crystals~\cite{FSS.2007,ST12}. Therefore, we believe that the Nakajima monomials appearing in a realization of general KR crystals will give a connection between the work of~\cite{HL16,KN15,Nakashima05,Nakashima14,RSW16}.
In particular, a key aspect of our construction is the Dynkin quiver being an oriented cycle as otherwise we obtain the crystal of the level-zero representation. In~\cite[Thm.~1.1(2)]{RSW16}, there is the condition that the Dynkin quiver have no oriented cycles, and so our results give evidence there are interesting cluster variables for the corresponding cluster algebra with an oriented cycle.
Furthermore, our construction gives a further connection between the Nakajima monomial model, the abacus model, ribbon partition model, (multi)partition model, and quiver varieties that was explored in~\cite{ST14,Tingley08}. Indeed, this allows us to construct the rank-level duality on Nakajima monomials and it suggests a crystal structure on cyclic quiver varieties~\cite{Nakajima01II}.


We now give one potential application of our results. To do so, recall geometric crystals introduced by Berenstein and Kazhdan~\cite{BK00,BK07}, where the $U_q(\fsl_n)$-crystal structure of $B(\lambda)$ is lifted to rational actions on algebraic varieties. This was generalized to a lifting of certain $U_q(\asl_n)$-crystals using Schubert varieties~\cite{Nakashima05}. Nakashima has lifted Nakajima monomials to describe the decoration function of geometric crystals (specifically as generalized minors) in~\cite{Nakashima14} and made a connection to the polyhedral model in~\cite{Nakashima13}. There is also a lifting of the $U_q'(\asl_n)$-crystal $B_{\infty}$ to the geometric setting given in~\cite{KNO08} and the coherent limits of $\{ B^{r,s} \}_{s=1}^{\infty}$ in~\cite{MN16}. The geometric $R$-matrix has also been studied~\cite{KNO10,Yamada01} and has been used to relate a quotient of the liftings of $B_{\infty}$ to the unipotent loop group in~\cite{LP11}. Recently, Frieden has given a lifting of $B^{r,s}$ for $\asl_n$ to the Grassmannian $\operatorname{Gr}(n-r,n)$~\cite{Frieden17} and the corresponding combinatorial $R$-matrix~\cite{Frieden17II}.

We expect that the natural extension of our results to the coherent limit can be lifted to a statement on geometric crystals, connecting the results of Nakashima with the work on geometric analogs of $B_{\infty}$, see Remark~\ref{remark:infinity_case}. Furthermore, we believe our results could be used to construct a geometric lifting of the path model embeddings and give a geometric lifting of highest weight $U_q(\asl_n)$-crystals. We also believe that our results could give a connection to the cluster algebra (geometric) $R$-matrices that were recently introduced in~\cite{ILP16}.

This paper is organized as follows. In Section~\ref{sec:background}, we give background on crystals, KR crystals, Nakajima monomials, and the Kyoto path model.
In Section~\ref{sec:monomial_model}, we construct a model for $B^{1,s}$ using Nakajima monomials.
In Section~\ref{sec:Kyoto_path}, we give a method to construct tensor products of KR crystals as a map on Nakajima monomials and relate our construction to the Kyoto path model.
In Section~\ref{sec:extensions}, we describe the relationship between our model and other models for highest weight $U_q(\asl_n)$-crystals, generalizations of our model to $B^{r,1}$, and (possible) extensions to other types.

\subsection*{Acknowledgements}

The authors would like to thank Peter Tingley, Rinat Kedem, and Bolor Turmunkh for valuable discussions. The authors would like to thank Masato Okado, Ben Salisbury, and Anne Schilling for comments on earlier drafts of this paper. The authors thank the anonymous referee for many useful comments and improvements to this paper. This work benefited from computations using \textsc{SageMath}~\cite{combinat,sage}.

The majority of this work was done while the authors were at the University of Minnesota.

\section{Background}
\label{sec:background}


\subsection{Crystals}

\begin{figure}
\begin{center}
\begin{tikzpicture}[scale=0.7]
\draw (6.7 cm, -1 cm) -- (3 cm,-1 cm) -- (1 cm,0) -- (7.5 cm, 1 cm) -- (14 cm, 0) -- (12 cm, -1 cm) -- (8.3 cm, -1 cm);
\draw[style=dashed] (6.7 cm,-1 cm) -- (8.3 cm,-1 cm);
\draw[fill=white] (7.5 cm, 1 cm) circle (.25cm) node[above=4pt]{$0$};
\draw[fill=white] (1 cm, 0 cm) circle (.25cm) node[left=4pt]{$1$};
\draw[fill=white] (3 cm, -1 cm) circle (.25cm) node[below=4pt]{$2$};
\draw[fill=white] (6 cm, -1 cm) circle (.25cm) node[below=4pt]{$3$};
\draw[fill=white] (9 cm, -1 cm) circle (.25cm) node[below=4pt]{$n-3$};
\draw[fill=white] (12 cm, -1 cm) circle (.25cm) node[below=4pt]{$n-2$};
\draw[fill=white] (14 cm, 0 cm) circle (.25cm) node[right=4pt]{$n-1$};
\end{tikzpicture}
\end{center}
\caption{The Dynkin diagram of $\asl_n$.}
\label{fig:dynkin_diagram}
\end{figure}

Let $\asl_n$ be the affine Kac--Moody Lie algebra of type $A_{n-1}^{(1)}$ with index set $I = \ZZ/n\ZZ = \{0, 1, \dotsc, n-1\}$, Cartan matrix $(a_{ij})_{i,j \in I}$, simple roots $\{\alpha_i\}_{i \in I}$, simple coroots $\{h_i\}_{i \in I}$, fundamental weights $\{\Lambda_i\}_{i \in I}$, weight lattice $P := \operatorname{span}_{\ZZ} \{\Lambda_i \mid i \in I\}$, dual weight lattice $P^{\vee}$, canonical pairing $\langle\ ,\ \rangle \colon P^{\vee} \times P \to \ZZ$ given by $\inner{h_i}{\alpha_j} = a_{ij}$, and quantum group $U_q(\asl_n)$.
See Figure~\ref{fig:dynkin_diagram} for the Dynkin diagram of $\asl_n$.
Let $P^+ := \operatorname{span}_{\ZZ_{\geq 0}} \{\Lambda_i \mid i \in I\}$ denote the dominant integral weights.
Note that $\fsl_n$ is the canonical simple Lie algebra given by the index set $I_0 := I \setminus \{0\}$. Let $\clfw_i$ denote the natural projection of $\Lambda_i$ onto the weight lattice $\cl{P}$ of $\fsl_n$.

Let $c = h_0 + h_1 + \cdots + h_{n-1}$ denote the canonical central element of $\asl_n$. We define the \defn{level} of a weight $\lambda$ as $\inner{c}{\lambda}$. Let $P_s^+ := \{ \lambda \in P^+ \mid \inner{c}{\lambda} = s \}$ denote the set of level $s$ weights.

We write $U_q'(\asl_n) := U_q([\asl_n, \asl_n])$, and let $\delta := \alpha_0 + \alpha_1 + \cdots + \alpha_{n-1}$ denote the null root. Note that the $U_q'(\asl_n)$ fundamental weights and simple roots are also given by $\{\Lambda_i\}_{i \in I}$ and $\{\alpha_i\}_{i \in I}$, respectively, but are considered in the weight lattice $P / \ZZ \delta$.

An \defn{abstract $U_q(\asl_n)$-crystal} is a set $B$ with \defn{crystal operators} $e_i, f_i \colon B \to B \sqcup \{0\}$ for $i \in I$, statistics $\varepsilon_i, \varphi_i \colon B \to \ZZ \sqcup \{-\infty\}$ for $i \in I$, and \defn{weight function} $\wt \colon B \to P$  that satisfy the following conditions for all $i \in I$:
\begin{enumerate}
\item[(1)] $\varphi_i(b) = \varepsilon_i(b) + \inner{h_i}{\wt(b)}$ for all $b \in B$ and $i \in I$;
\item[(2)] if $e_i b \neq 0$ for $b \in B$, then
\begin{enumerate}
  \item $\varepsilon_i(e_i b) = \varepsilon_i(b) - 1$,
  \item $\varphi_i(e_i b) = \varphi_i(b) + 1$,
  \item $\wt(e_i b) = \wt(b) + \alpha_i$;
\end{enumerate}
\item[(3)] if $f_i b \neq 0$ for $b \in B$, then
\begin{enumerate}
  \item $\varepsilon_i(f_i b) = \varepsilon_i(b) + 1$,
  \item $\varphi_i(f_i b) = \varphi_i(b) - 1$,
  \item $\wt(f_i b) = \wt(b) - \alpha_i$;
\end{enumerate}
\item[(4)] $f_i b = b'$ if and only if $b = e_i b'$ for $b, b' \in B$ and $i \in I$;
\item[(5)] if $\varphi_i(b) = -\infty$ for $b \in B$, then $e_i b = f_i b = 0$.
\end{enumerate}
An abstract $U_q'(\asl_n)$-crystal is a $U_q(\asl_n)$-crystal except the weight function takes values in $P/\ZZ\delta$ rather than $P$.

Define
\[
\varepsilon(b) = \sum_{i \in I} \varepsilon_i(b) \Lambda_i,
\qquad\qquad
\varphi(b) = \sum_{i \in I} \varphi_i(b) \Lambda_i.
\]
We say an element $b \in B$ is \defn{highest weight} if $e_i b = 0$ for all $i \in I$. If $B$ is an abstract $U_q'(\asl_n)$-crystal, then we say $b \in B$ is \defn{classically highest weight} if $e_i b = 0$ for all $i \in I_0$.
We say an abstract $U_q(\asl_n)$-crystal is \defn{regular} if
\[
\varepsilon_i(b) = \max \{ k \mid e_i^k b \neq 0 \},
\qquad \qquad \varphi_i(b) = \max \{ k \mid f_i^k b \neq 0 \}.
\]
\begin{remark}
The term regular is sometimes called \defn{seminormal} in the literature.
\end{remark}

We call an abstract $U_q(\asl_n)$-crystal $B$ a \defn{$U_q(\asl_n)$-crystal} if $B$ is the crystal basis of some $U_q(\asl_n)$-module. Similarly for $U_q'(\asl_n)$.

Kashiwara showed in~\cite{K91} that the irreducible highest weight $U_q(\g)$-module $V(\lambda)$ admits a crystal basis, where $\g$ is a symmetrizable Kac--Moody Lie algebra and $\lambda$ is a dominant integral weight. We denote this crystal basis by $B(\lambda)$, and let $u_{\lambda} \in B(\lambda)$ denote the unique highest weight element, which is the unique element of weight $\lambda$.

We define the \defn{tensor product} of abstract $U_q(\asl_n)$-crystals $B_1$ and $B_2$ as the crystal $B_2 \otimes B_1$ that is the Cartesian product $B_2 \times B_1$ with the crystal structure
\begin{align*}
e_i(b_2 \otimes b_1) & = \begin{cases}
e_i b_2 \otimes b_1 & \text{if } \varepsilon_i(b_2) > \varphi_i(b_1), \\
b_2 \otimes e_i b_1 & \text{if } \varepsilon_i(b_2) \leq \varphi_i(b_1),
\end{cases}
\\ f_i(b_2 \otimes b_1) & = \begin{cases}
f_i b_2 \otimes b_1 & \text{if } \varepsilon_i(b_2) \geq \varphi_i(b_1), \\
b_2 \otimes f_i b_1 & \text{if } \varepsilon_i(b_2) < \varphi_i(b_1),
\end{cases}
\\ \varepsilon_i(b_2 \otimes b_1) & = \max(\varepsilon_i(b_1), \varepsilon_i(b_2) - \inner{h_i}{\wt(b_1)}),
\\ \varphi_i(b_2 \otimes b_1) & = \max(\varphi_i(b_2), \varphi_i(b_1) + \inner{h_i}{\wt(b_2)}),
\\ \wt(b_2 \otimes b_1) & = \wt(b_2) + \wt(b_1).
\end{align*}
\begin{remark}
Our tensor product convention is opposite of Kashiwara's~\cite{K91}.
\end{remark}

Let $B_1$ and $B_2$ be two abstract $U_q(\g)$-crystals. A \defn{crystal morphism} $\psi \colon B_1 \to B_2$ is a map $B_1 \sqcup \{0\} \to B_2 \sqcup \{0\}$ with $\psi(0) = 0$ such that the following properties hold for all $b \in B_1$:
\begin{itemize}
\item[(1)] If $\psi(b) \in B_2$, then $\wt\bigl(\psi(b)\bigr) = \wt(b)$, $\varepsilon_i\bigl(\psi(b)\bigr) = \varepsilon_i(b)$, and $\varphi_i\bigl(\psi(b)\bigr) = \varphi_i(b)$.
\item[(2)] We have $\psi(e_i b) = e_i \psi(b)$ if $\psi(e_i b) \neq 0$ and $e_i \psi(b) \neq 0$.
\item[(3)] We have $\psi(f_i b) = f_i \psi(b)$ if $\psi(f_i b) \neq 0$ and $f_i \psi(b) \neq 0$.
\end{itemize}
An \defn{embedding} (resp.\ \defn{isomorphism}) is a crystal morphism such that the induced map $B_1 \sqcup \{0\} \to B_2 \sqcup \{0\}$ is an embedding (resp.\ bijection). A crystal morphism is \defn{strict} if it commutes with all crystal operators. Note that a crystal isomorphism is strict.

\subsection{Nakajima monomials}

Next, we recall the Nakajima monomial realization of crystals following~\cite{ST14}.\footnote{There is a minor typo in~\cite{ST14} with the partial sums in $k_e$ being over $s > k$ instead of $s \geq k$.}

Let $\mcM$ denote the set of Laurent monomials in the commuting variables $\{Y_{i,k}\}_{i \in I, k \in \ZZ}$. Fix an integer $K$, and then fix integers $c_{i,i+1}$ and $c_{i+1,i}$ such that $K = c_{i,i+1} + c_{i+1,i}$ for all $i \in I$, where all indices are taken mod $n$. For a monomial $m = \prod_{i \in I} \prod_{k \in \ZZ} Y_{i,k}^{y_{i,k}}$, define
\begin{align*}
\varepsilon_i(m) & = -\min_{k \in \ZZ} \sum_{s \geq k} y_{i,s},
& k_e(m) & = \max \left\{ k \exmid \varepsilon_i(m) = - \sum_{s \geq k} y_{i,s} \right\},
\\ \varphi_i(m) & = \max_{k \in \ZZ} \sum_{s \leq k} y_{i,s},
& k_f(m) & = \min \left\{ k \exmid \varphi_i(m) = \sum_{s \leq k} y_{i,s} \right\},
\end{align*}
\[
\wt(m) = \sum_{\substack{i \in I \\ k \in \ZZ}} y_{i,k} \Lambda_i.
\]
Note that $k_e$ (resp.~$k_f$) is undefined if $\varepsilon_i(m) = 0$ (resp.~$\varphi_i(m) = 0$).
Define the crystal operators $e_i,f_i \colon \mcM \to \mcM \sqcup \{ 0 \}$ by
\begin{align*}
e_i(m) & = \begin{cases} 0 & \text{if } \varepsilon_i(m) = 0, \\ m A_{i, k_e(m) - K} & \text{if } \varepsilon_i(m) > 0, \end{cases}
& f_i(m) & = \begin{cases} 0 & \text{if } \varphi_i(m) = 0, \\ m A_{i, k_f(m)}^{-1} & \text{if } \varphi_i(m) > 0, \end{cases}
\end{align*}
where
\[
A_{i,k} = Y_{i,k} Y_{i,k+K} Y_{i-1,k+c_{i-1,i}}^{-1} Y_{i+1,k+c_{i+1,i}}^{-1}.
\]

\begin{figure}
\[
\begin{tikzpicture}[>=latex,xscale=3,yscale=1.5]
\node (hw) at (0,0) {$Y_{1,0} Y_{2,0}$};
\node (1) at (1,1) {$Y_{1,1}^{-1} Y_{2,0}^2$};
\node (21) at (2,1) {$Y_{2,0} Y_{2,1}^{-1}$};
\node (221) at (3,1) {$Y_{1,1} Y_{2,1}^{-2}$};
\node (2) at (1,-1) {$Y_{1,0} Y_{1,1} Y_{2,1}^{-1}$};
\node (12) at (2,-1) {$Y_{1,0} Y_{1,2}^{-1}$};
\node (112) at (3,-1) {$Y_{1,1}^{-1} Y_{1,2}^{-1} Y_{2,0}$};
\node (lw) at (4,0) {$Y_{1,2}^{-1} Y_{2,1}^{-1}$};
\draw[->,blue] (hw) -- node[midway,above left,black] {\small $1$} (1);
\draw[->,red] (hw) -- node[midway,below left,black] {\small $2$} (2);
\draw[->,red] (1) -- node[midway,above,black] {\small $2$} (21);
\draw[->,red] (21) -- node[midway,above,black] {\small $2$} (221);
\draw[->,blue] (2) -- node[midway,below,black] {\small $1$} (12);
\draw[->,blue] (12) -- node[midway,below,black] {\small $1$} (112);
\draw[->,blue] (221) -- node[midway,above right,black] {\small $1$} (lw);
\draw[->,red] (112) -- node[midway,below right,black] {\small $2$} (lw);
\end{tikzpicture}
\]
\caption{The crystal $\mcM(\clfw_1 + \clfw_2)$ in type $A_2$.}
\label{fig:monomial_example}
\end{figure}

We note that the crystal structure of Kashiwara~\cite{K03II} is when $K = 1$, and that of Nakajima~\cite{Nakajima03} is when $c_{i,i+1} = c_{i+1,i} = 1$ for all $i \in I$ (so $K = 2$). Note that in the case of Nakajima, odd length cycles are prohibited in the Dynkin diagram; that is, we can only consider types $\asl_{2k}$. 

Let $\mcM(m)$ denote the closure of $m$ under the crystal operators $e_i$ and $f_i$ for all $i \in I$.

\begin{theorem}[{Sam--Tingley~\cite{ST14}}]
\label{thm:highest_weight_monomials}
For $\lambda \in P^+$ and any $(k_i)_{i \in I}$, we have 
\[
\mcM\left( \prod_{i \in I} Y_{i,k_i}^{\inner{h_i}{\lambda}} \right) \iso B(\lambda).
\]
\end{theorem}

See Figure~\ref{fig:monomial_example} for an example that is isomorphic to $B(\clfw_1 + \clfw_2)$ in type $A_2$.

Note that for a monomial $m$ given in Theorem~\ref{thm:highest_weight_monomials}, the resulting crystal $\mcM(m)$ is a regular crystal.
For convenience, we define
\[
Y_{\lambda} := \prod_{i \in I} Y_{i, 0}^{\inner{h_i}{\lambda}},
\]
and we denote $\mathbf{1} = Y_0$, where $0 \in P^+$. We also define $\mcM(\lambda) := \mcM(Y_{\lambda})$.

%

There is also another set of variables
\begin{equation}
\label{eq:X_variables}
X_{i,k} := Y_{i-1,k+1}^{-1} Y_{i,k}
\end{equation}
that was used in the description of $\mcM(\lambda)$ in~\cite[Cor.~4.9]{Kim05} and $\mcM(\infty)$ in~\cite[Cor.~4.3]{KKS07} for type $A_n^{(1)}$.

\begin{lemma}
\label{lemma:positivity}
Let $\lambda \in P^+$ and $c_{ij} \in \ZZ_{\geq 0}$. The monomials in $\mcM(\lambda)$ only contain variables $Y_{i,k}$ with $k \geq 0$ and the exponent of $Y_{i,0}$ is nonnegative.
\end{lemma}

\begin{proof}
Recall that $\mcM(\lambda)$ is generated by applying $f_i$, for $i \in I$, to $Y_{\lambda}$.
Note that $A_{i,k}$ for $k \geq 0$ does not contain $y_{i',k'}$ for $k' < 0$.
We proceed using an induction on depth.
Note for $\lambda \in P^+$, since $\varphi_i(m) > 0$, by our induction assumption we must have $k_f(m) \geq 0$.
Furthermore, $\varphi_i(m) > 0$ and $k_f(m) = 0$ imply there exists a $Y_{i,0}$ occurring with a positive exponent.
Hence, the claim follows by the definition of the crystal operators.
(See also~\cite[Cor.~4.9]{Kim05}.) 
\qed\end{proof}


Unless otherwise stated, we consider $c_{ij} = 1$ if $(i, j) = (n, 0)$ or $i < j$ when $(i, j) \neq (0, n)$ and $c_{ij} = 0$ otherwise (hence $K = 1$). Note that this corresponds to orienting the Dynkin diagram into an ordered cycle, where we draw an arrow $i \to i+1$ implying $c_{i,i+1} = 1$ and the other values $c_{i,j}$ for $j \neq i \pm 1$ do not affect the crystal structure.
Furthermore, $K=1$ means our $e_i$ crystal operators can alternatively be given by
\[
n_e(m) = \max \left\{ k \exmid \varepsilon_i(m) = - \sum_{s > k} y_{i,s} \right\},
\qquad\qquad
e^!_i(m) = \begin{cases} 0 & \text{if } \varepsilon_i(m) = 0, \\ m A_{i, n_e(m)} & \text{if } \varepsilon_i(m) > 0. \end{cases}
\]
By using the relation $\varepsilon_i(m) = \varphi_i(m) - \inner{h_i}{\wt(m)}$, we can rewrite
\begin{equation}
\label{eq:alt_ke}
n_e(m) = \max \left\{ k \exmid \varphi_i(M) = \sum_{s \leq k} y_{is} \right\}.
\end{equation}
(See also, {\it e.g.},~\cite{K03II}.) Therefore, our crystal structure is equivalent to that given in~\cite{KKS07}.

\subsection{Kirillov--Reshetikhin crystals}
\label{sec:KR_crystals}

A \defn{Kirillov--Reshetikhin (KR) module} $W^{r,s}$, where $r \in I_0$ and $s \in \ZZ_{> 0}$, is a particular irreducible finite-dimensional $U_q'(\asl_n)$-module that has many remarkable properties. KR modules are classified by their Drinfel'd polynomials, and $W^{r,s}$ is the minimal affinization of the highest weight $U_q(\fsl_n)$-representation $V(s \clfw_r)$~\cite{CP95,CP98}. In particular, it was shown in~\cite{KKMMNN92} that the KR module $W^{r,s}$ admits a crystal basis $B^{r,s}$ called a \defn{Kirillov--Reshetikhin (KR) crystal}.
Another property is that the KR crystal $B^{r,s}$ is a \defn{perfect crystal of level $s$}, which means it satisfies the following conditions:
\begin{enumerate}
\item $B^{r,s} \otimes B^{r,s}$ is connected.
\item $\cl{\wt}(b) \in s \clfw_r + \sum_{i \in I_0} \ZZ_{\leq 0} \alpha_i$ for all $b \in B^{r,s}$.
\item $\inner{c}{\varepsilon(b)} \geq s$ for all $b \in B^{r,s}$.
\item For all $\lambda \in P_s^+$, there exist unique elements $b_{\lambda}, b^{\lambda} \in B^{r,s}$ such that
\[
\varepsilon(b_{\lambda}) = \lambda = \varphi(b^{\lambda}).
\]
\end{enumerate}

Additionally, we have $B^{r,s} \iso B(s\clfw_r)$ as $U_q(\fsl_n)$-crystals, and so $B^{r,s}$ can be modeled by semistandard tableaux of shape $(s^r)$. Next, recall that the Dynkin diagram automorphism $i \mapsto i+1 \mod n$ induces a (twisted) $U_q(\asl_n)$-crystal isomorphism $\pr \colon B^{r,s} \to B^{r,s}$ called the \defn{promotion isomorphism}~\cite{Shimozono02}. On semistandard tableaux~\cite{KN94}, the map $\pr$ is the (weak\footnote{See, \textit{e.g.},~\cite{BST10,Kus13} for a definition.}) promotion operator of Sch\"utzenberger~\cite{Sch72}. Hence, $B^{r,s}$ is a regular crystal, and we define the remaining crystal structure on $B^{r,s}$ by
\begin{align*}
e_0 & = \pr^{-1} \circ \; e_1 \circ \pr,
\\ f_0 & = \pr^{-1} \circ \; f_1 \circ \pr,
\\ \wt(b) & = \cl{\wt}(b) + k_0 \Lambda_0,
\end{align*}
where $k_0$ is such that $\inner{\wt(b)}{c} = 0$ ({\it i.e.}, it is a level 0 weight).

We will be focusing on the KR crystal $B^{1,s}$, which is the crystal corresponding to $W^{1,s}$ and can be modeled using the \defn{vector realization}. Indeed, we have
\[
B^{1,s} = \left\{ (x_1, \dotsc, x_n) \exmid x_1, \dotsc, x_n \in \ZZ_{\geq 0},\, \sum_{i=1}^n x_i = s \right\}
\]
with the crystal structure
\begin{subequations}
\label{eq:vector_crystal}
\begin{align}
e_i(x_1, \dotsc, x_n) & = \begin{cases} 0 & \text{if } x_{i+1} = 0, \\ (x_1, \dotsc, x_i + 1, x_{i+1} - 1, \dotsc, x_n) & \text{if } x_{i+1} > 0, \end{cases}
\\ f_i(x_1, \dotsc, x_n) & = \begin{cases} 0 & \text{if } x_i = 0, \\ (x_1, \dotsc, x_i - 1, x_{i+1} + 1, \dotsc, x_n) & \text{if } x_i > 0, \end{cases}
\\ \varepsilon_i(x_1, \dotsc, x_n) & = x_{i+1},
\\ \varphi_i(x_1, \dotsc, x_n) & = x_i,
\\ \wt(x_1, \dotsc, x_n) & = \sum_{i \in I} (x_i - x_{i+1}) \Lambda_i,
\end{align}
\end{subequations}
where all indices are understood mod $n$. Note that $B^{1,s}$ is a regular crystal and is naturally identified with semistandard tableaux of shape $(s)$, where $x_i$ equals the number of $i$'s in a tableau.

We will also need the \defn{affinization} of a $U_q'(\asl_n)$-crystal $B$, which is defined as follows.
The affinization of $B$ is the $U_q(\asl_n)$-crystal $\widehat{B} =\{ b(k) \mid b \in B, k \in \ZZ \}$, whose crystal structure is given by
\begin{align*}
e_i\bigl(b(k)\bigr) & = \begin{cases} (e_0 b)(k+1) & \text{if } i = 0, \\ (e_i b)(k) & \text{if } i \neq 0, \end{cases}
\\ f_i\bigl(b(k)\bigr) & = \begin{cases} (f_0 b)(k-1) & \text{if } i = 0, \\ (f_i b)(k) & \text{if } i \neq 0, \end{cases}
\\ \varepsilon_i\bigl(b(k)\bigr) & = \varepsilon_i(b),
\\ \varphi_i\bigl(b(k)\bigr) & = \varphi_i(b),
\\ \wt\bigl(b(k)\bigr) & = \wt(b) + k \delta.
\end{align*}


\subsection{Kyoto path model}

We recall some of the results of~\cite{KKMMNN91,KKMMNN92,OSS03IV}, which give a model for highest weight $U_q(\asl_n)$-crystals using KR crystals.

\begin{theorem}
\label{thm:kyoto_embedding}
Let $\lambda \in P^+$ be a level $s$ weight. Let $B$ be a perfect crystal of level $s$. Let $b^{\lambda} \in B$ be the unique element such that $\varphi(b^{\lambda}) = \lambda$. Let $\mu = \varepsilon(b^{\lambda})$. The morphism
\[
\widehat{\Psi} \colon B(\lambda) \to \widehat{B} \otimes B(\mu)
\]
defined by $u_{\lambda} \mapsto b^{\lambda}(0) \otimes u_{\mu}$ is a strict $U_q(\asl_n)$-crystal embedding.
\end{theorem}

By iterating $\widehat{\Psi}$, we obtain the \defn{Kyoto path model}. For a level $s$ weight $\lambda$, we can construct a model for $B(\lambda)$ by
\[
\widehat{\Psi}^{(+\infty)} \colon B(\lambda) \to \widehat{B}^{1,s} \otimes \widehat{B}^{1,s} \otimes \cdots
\]
since $B^{1,s}$ is a perfect crystal of level $s$.
Note that we consider $e_i b = 0$ if the action would be otherwise undefined in $\widehat{B}^{1,s} \otimes \widehat{B}^{1,s} \otimes \cdots$.
Furthermore, $\widehat{\Psi}^{(+\infty)}(u_{\lambda})$ is eventually cyclic and, for any $b \in B(\lambda)$, the element $\widehat{\Psi}^{(+\infty)}(b)$ only differs from $\widehat{\Psi}^{(+\infty)}(u_{\lambda})$ in a finite number of factors. Therefore, for any element $b$, we can take Theorem~\ref{thm:kyoto_embedding} iterated $N \gg 1$ times (which depends on $b$) $\widehat{\Psi}^{(N)}(b)$ to define the crystal structure on the Kyoto path model using only the KR crystal $B^{1,s}$.

We note that there are analogous results for $U_q'(\asl_n)$-crystals by considering the branching rule from $U_q(\asl_n)$ to $U_q'(\asl_n)$. In particular, there exists a $U_q'(\asl_n)$-crystal isomorphism
\begin{equation}
\label{eq:kyoto_prime}
\Psi \colon B(\lambda) \to B \otimes B(\mu)
\end{equation}
defined by $u_{\lambda} \mapsto b^{\lambda} \otimes u_{\mu}$.

%

\section{Monomial realization of \texorpdfstring{$B^{1,s}$}{B1s}}
\label{sec:monomial_model}

In this section, we describe the construction of $B^{1,s}$ using Kashiwara's crystal structure~\cite{K03II} of Nakajima monomials.

Define
\[
\mcM^{1,s} := \left\{ \prod_{i=1}^n Y_{i-1,1}^{-x_i} Y_{i,0}^{x_i} \exmid x_1, \dotsc, x_n \in \ZZ_{\geq 0}, \sum_{i=1}^n x_i = s \right\}.
\]

\begin{theorem}
\label{thm:single_row_isomorphism}
We have
\[
\mcM^{1,s} \iso B^{1,s}
\]
as $U_q'(\asl_n)$-crystals.
\end{theorem}

\begin{proof}
Let $\Phi \colon B^{1,s} \to \mcM^{1,s}$ be the map
\[
(x_1, \dotsc, x_n) \mapsto \prod_{i=1}^n Y_{i-1,1}^{-x_i} Y_{i,0}^{x_i},
\]
and it is clear that $\Phi$ is a weight preserving bijection. Thus it remains to show that $\Phi$ commutes with the crystal operators since $B^{1,s}$ and $\mcM^{1,s}$ are regular crystals.

We restrict to the variables only containing $i$ since these completely determine the action of $e_i$ and $f_i$. It is sufficient to consider
\[
m = Y_{i,1}^{k_1}\, Y_{i,0}^{k_2}
\]
with $k_1\leq 0 \leq k_2$. Then
\begin{subequations}
\label{eq:ep_phi_M1s}
\begin{align}
\varphi_i(m) & = \max\left\{\sum_{s\leq k} y_{i,s} \exmid k \in \ZZ \right\} = \max\{k_2, k_2+k_1\} = k_2 \geq 0,
\\
\varepsilon_i(m) & = -\min\left\{\sum_{s \geq k} y_{i,s} \exmid k \in \ZZ \right\} = \max\{-k_2-k_1,-k_1\} = -k_1 \geq 0,
\end{align}
\end{subequations}
since $k_1 \leq 0 \leq k_2$. Note that $\varphi_i(b) = 0$ (so $f_i(b) = 0$) if and only if $f_i \Phi(b) = 0$, and similarly for $\varepsilon_i$ and $e_i$.

Hence 
\begin{align*}
k_f(m) & = \min\left\{k \exmid \varphi_i(m) = \sum_{s\leq k} y_{i,s}\right\}
= \min\left\{k \exmid k_2 = \sum_{s \leq k} y_{i,s}\right\} = 0,
\\ k_e(m) & = \max\left\{k \exmid -\varepsilon_i(m) = \sum_{s \geq k} y_{i,s}\right\}
= \max\left\{k \exmid k_1 = \sum_{s \geq k} y_{i,s}\right\} = 1.
\end{align*}
Therefore, we have
\begin{align*}
A_{i,k_f(m)} = A_{i,k_e(m)-1} = A_{i,0} &= 
Y_{i,0} \, Y_{i,1} \, Y_{i-1,c_{i-1,i}}^{-1} \, Y_{i+1,c_{i+1,i}}^{-1} 
\\ &=Y_{i,0} \, Y_{i,1} \, Y_{i-1,1}^{-1} \, Y_{i+1,0}^{-1}.
\end{align*}
Since $\varphi_i(m) > 0$ and $\varepsilon_i(m) > 0$, we have
\begin{align*}
f_i(m) & =
A_{i,k_f(m)}^{-1} m =
Y_{i,0}^{k_2-1} \, Y_{i,1}^{k_1-1} \, Y_{i-1,1}^{} \, Y_{i+1,0}^{},
\\ e_i(m) & =
A_{i,k_e(m)-1} m =
Y_{i,0}^{k_2+1} \, Y_{i,1}^{k_1+1} \, Y_{i-1,1}^{-1} \, Y_{i+1,0}^{-1}. 
\end{align*}

For $i \in I$, we have
\begin{align*}
f_i \bigl( \Phi (x_1, \dotsc, x_n) \bigr) &=  f_i \left( \prod_{k=1}^n  Y_{k-1,1}^{-x_k} Y_{k,0}^{x_k} \right)
\\ &= f_i\left(Y_{i-1,1}^{-x_i} Y_{i,0}^{x_i}   \,\, Y_{i,1}^{-x_{i+1}} Y_{i+1,0}^{x_{i+1}}  \right) \prod_{k\neq i,i+1}  Y_{k-1,1}^{-x_k} Y_{k,0}^{x_k} 
\\ &= f_i\left(Y_{i,1}^{-x_{i+1}} Y_{i,0}^{x_i} \right)  \,\, Y_{i-1,1}^{-x_i} Y_{i+1,0}^{x_{i+1}}  \prod_{k\neq i,i+1}  Y_{k-1,1}^{-x_k} Y_{k,0}^{x_k} 
\\ &= \left(Y_{i,1}^{-x_{i+1}-1} Y_{i,0}^{x_i-1}   \,\,\,\, Y_{i-1,1} Y_{i+1,0} \right) 
 \,\, Y_{i-1,1}^{-x_i} Y_{i+1,0}^{x_{i+1}}  \prod_{k\neq i,i+1}  Y_{k-1,1}^{-x_k} Y_{k,0}^{x_k} 
\\ &= Y_{i-1,1}^{-x_i+1} \, Y_{i,0}^{x_i-1} \, Y_{i,1}^{-x_{i+1}-1} \, Y_{i+1,0}^{x_{i+1}+1}  \prod_{k\neq i,i+1}  Y_{k-1,1}^{-x_k} Y_{k,0}^{x_k} 
\\ &= \Phi(\dotsc, x_i - 1, x_{i+1} + 1, \ldots)
\\ &= \Phi\bigl( f_i (x_1, \dotsc, x_n) \bigr)
\end{align*}
where all indices are taken mod $n$. Similarly, we have $e_i\bigl( \Phi (x_1, \dotsc, x_n) \bigr) = \Phi\bigl( e_i (x_1, \dotsc, x_n) \bigr)$.
\qed\end{proof}

Note that we can define
\begin{equation}
\label{eq:image_X}
\Phi(x_1, \dotsc, x_n) = X_{1,0}^{x_1} X_{2,0}^{x_2} \dotsm X_{n,0}^{x_n},
\end{equation}
where the variables $X_{i,k}$ are given by Equation~\eqref{eq:X_variables}.

\begin{example}
The crystal $\mcM^{1,3}$ for $\asl_3$ is given by Figure~\ref{fig:asl3_M13_ex}.
\end{example}

\begin{figure}
\[
\begin{tikzpicture}[>=latex,line join=bevel,xscale=1.0, yscale=0.7, every node/.style={scale=0.9}]
\node (node_9) at (193.0bp,314.0bp) [draw,draw=none] {$Y_{0,1}^{-1} Y_{1,0} Y_{1,1}^{-2} Y_{2,0}^{2} $};
  \node (node_8) at (38.0bp,314.0bp) [draw,draw=none] {$Y_{0,0} Y_{0,1}^{-2} Y_{1,0}^{2} Y_{2,1}^{-1} $};
  \node (node_7) at (258.0bp,238.0bp) [draw,draw=none] {$Y_{1,1}^{-3} Y_{2,0}^{3} $};
  \node (node_6) at (115.0bp,86.0bp) [draw,draw=none] {$Y_{0,0}^{2} Y_{1,1}^{-1} Y_{2,0} Y_{2,1}^{-2} $};
  \node (node_5) at (113.0bp,390.0bp) [draw,draw=none] {$Y_{0,1}^{-2} Y_{1,0}^{2} Y_{1,1}^{-1} Y_{2,0} $};
  \node (node_4) at (65.0bp,466.0bp) [draw,draw=none] {$Y_{0,1}^{-3} Y_{1,0}^{3} $};
  \node (node_3) at (195.0bp,162.0bp) [draw,draw=none] {$Y_{0,0} Y_{1,1}^{-2} Y_{2,0}^{2} Y_{2,1}^{-1} $};
  \node (node_2) at (40.0bp,162.0bp) [draw,draw=none] {$Y_{0,0}^{2} Y_{0,1}^{-1} Y_{1,0} Y_{2,1}^{-2} $};
  \node (node_1) at (109.0bp,238.0bp) [draw,draw=none] {$Y_{0,0} Y_{0,1}^{-1} Y_{1,0} Y_{1,1}^{-1} Y_{2,0} Y_{2,1}^{-1} $};
  \node (node_0) at (67.0bp,10.0bp) [draw,draw=none] {$Y_{0,0}^{3} Y_{2,1}^{-3} $};
  \draw [black,<-] (node_4) ..controls (55.089bp,433.79bp) and (50.04bp,416.0bp)  .. (47.0bp,400.0bp) .. controls (41.735bp,372.3bp) and (39.328bp,338.9bp)  .. (node_8);
  \definecolor{strokecol}{rgb}{0.0,0.0,0.0};
  \pgfsetstrokecolor{strokecol}
  \draw (56.0bp,390.0bp) node {$0$};
  \draw [blue,->] (node_4) ..controls (78.345bp,444.43bp) and (92.105bp,423.21bp)  .. (node_5);
  \draw (104.0bp,428.0bp) node {$1$};
  \draw [blue,->] (node_5) ..controls (135.6bp,368.09bp) and (159.48bp,346.01bp)  .. (node_9);
  \draw (172.0bp,352.0bp) node {$1$};
  \draw [blue,->] (node_9) ..controls (211.27bp,292.2bp) and (230.41bp,270.41bp)  .. (node_7);
  \draw (242.0bp,276.0bp) node {$1$};
  \draw [red,->] (node_6) ..controls (101.66bp,64.427bp) and (87.895bp,43.214bp)  .. (node_0);
  \draw (106.0bp,48.0bp) node {$2$};
  \draw [red,->] (node_5) ..controls (95.71bp,375.2bp) and (87.71bp,368.45bp)  .. (81.0bp,362.0bp) .. controls (70.864bp,352.26bp) and (60.194bp,340.61bp)  .. (node_8);
  \draw (90.0bp,352.0bp) node {$2$};
  \draw [black,<-] (node_2) ..controls (42.331bp,124.78bp) and (44.723bp,98.505bp)  .. (49.0bp,76.0bp) .. controls (52.892bp,55.524bp) and (60.074bp,32.097bp)  .. (node_0);
  \draw (58.0bp,86.0bp) node {$0$};
  \draw [blue,->] (node_8) ..controls (54.256bp,294.63bp) and (68.344bp,278.99bp)  .. (81.0bp,266.0bp) .. controls (84.637bp,262.27bp) and (88.626bp,258.33bp)  .. (node_1);
  \draw (90.0bp,276.0bp) node {$1$};
  \draw [red,->] (node_1) ..controls (95.201bp,222.82bp) and (88.701bp,216.07bp)  .. (83.0bp,210.0bp) .. controls (73.496bp,199.88bp) and (62.971bp,188.38bp)  .. (node_2);
  \draw (92.0bp,200.0bp) node {$2$};
  \draw [black,<-] (node_8) ..controls (28.42bp,277.02bp) and (23.567bp,250.85bp)  .. (26.0bp,228.0bp) .. controls (28.175bp,207.57bp) and (34.094bp,184.12bp)  .. (node_2);
  \draw (35.0bp,238.0bp) node {$0$};
  \draw [blue,->] (node_1) ..controls (133.5bp,215.92bp) and (159.72bp,193.36bp)  .. (node_3);
  \draw (172.0bp,200.0bp) node {$1$};
  \draw [black,<-] (node_5) ..controls (111.62bp,337.21bp) and (109.86bp,271.34bp)  .. (node_1);
  \draw (121.0bp,314.0bp) node {$0$};
  \draw [black,<-] (node_9) ..controls (193.69bp,261.21bp) and (194.57bp,195.34bp)  .. (node_3);
  \draw (204.0bp,238.0bp) node {$0$};
  \draw [red,->] (node_7) ..controls (240.39bp,216.32bp) and (222.08bp,194.81bp)  .. (node_3);
  \draw (242.0bp,200.0bp) node {$2$};
  \draw [blue,->] (node_2) ..controls (55.611bp,142.7bp) and (69.761bp,126.72bp)  .. (83.0bp,114.0bp) .. controls (87.089bp,110.07bp) and (91.657bp,106.03bp)  .. (node_6);
  \draw (92.0bp,124.0bp) node {$1$};
  \draw [black,<-] (node_1) ..controls (111.07bp,185.21bp) and (113.71bp,119.34bp)  .. (node_6);
  \draw (122.0bp,162.0bp) node {$0$};
  \draw [red,->] (node_3) ..controls (172.4bp,140.09bp) and (148.52bp,118.01bp)  .. (node_6);
  \draw (174.0bp,124.0bp) node {$2$};
  \draw [red,->] (node_9) ..controls (169.07bp,291.92bp) and (143.46bp,269.36bp)  .. (node_1);
  \draw (170.0bp,276.0bp) node {$2$};
\end{tikzpicture}
\]
\caption{The crystal $\mcM^{1,3}$ for $\asl_3$.}
\label{fig:asl3_M13_ex}
\end{figure}

\section{Relation to the Kyoto path model}
\label{sec:Kyoto_path}

To define $\mcM^{1,s}$, we considered the crystal generated from $Y_{0,1}^{-s} Y_{1,0}^s$. However, by shifting the monomials, we can construct an isomorphism with the tensor product. Indeed, let $\tau_j$ be the map given by $Y_{i,k} \mapsto Y_{i,k+j}$ for all $i \in I$ and $k \in \ZZ$. Let $\mcM \cdot \mcM' = \{ m \cdot m' \mid m \in \mcM, m' \in \mcM' \}$, where $\cdot$ denotes the usual multiplication in $\mcM$, and endow this set with the usual crystal operators on Nakajima monomials.

\begin{lemma}
\label{lemma:product_isomorphism}
Let $\mcM$ and $\mcM'$ be two crystals of Nakajima monomials such that there exists a $J$ such that for all monomials
\[
\prod_{\substack{i \in I \\ k \in \ZZ}} Y_{i,k}^{y_{i,k}} \in \mcM,
\qquad\qquad
\prod_{\substack{i \in I \\ k \in \ZZ}} Y_{i,k}^{y'_{i,k}} \in \mcM',
\]
and for all $i \in I$, we have
\begin{itemize}
\item $y_{i,k} = 0$ for all $k > J$ with $y_{i,J} \leq 0$,
\item $y'_{i,k} = 0$ for all $k < J$ with $y'_{i,J} \geq 0$,
\end{itemize}
and there do not exist distinct monomials of the form
\[
\overline{m} \prod_{i \in I} Y_{i,J}^{y_{i,J}}, \ \overline{m} \prod_{i \in I} Y_{i,J}^{\check{y}_{i,J}} \in \mcM,
\qquad\qquad
\overline{m}' \prod_{i \in I} Y_{i,J}^{y'_{i,J}}, \ \overline{m}' \prod_{i \in I} Y_{i,J}^{\check{y}'_{i,J}} \in \mcM',
\]
such that $y_{i,J} + y'_{i,J} = \check{y}_{i,J} + \check{y}'_{i,J}$ for all $i \in I$.
Then there exists a crystal isomorphism
\[
\psi \colon \mcM \otimes \mcM' \to \mcM \cdot \mcM'
\]
given by $\psi(m \otimes m') = m \cdot m'$.
\end{lemma}

\begin{proof}
The final hypothesis guarantees that we can separate any monomial in $\mcM \otimes \mcM'$ into unique factors from $\mcM$ and $\mcM'$.
Indeed, the only place where there could be interaction between an element in $\mcM$ and $\mcM'$ is with the variables $Y_{i,J}$, where our condition that there exists an $i \in I$ such that $y_{i,J} + y'_{i,J} \neq \check{y}_{i,J} + \check{y}'_{i,J}$ implies the injectivity of $\psi$.
The map $\psi$ is clearly surjective, and thus $\psi$ is a bijection. It remains to show that $\psi$ is a crystal morphism.

From our assumptions on $\mcM$ and $\mcM'$, we have
\begin{subequations}
\label{eq:products}
\begin{align}
\label{eq:product_epsilon}
\sum_{p \geq k} (y_{i,p} + y'_{i,p}) & =
\begin{cases}
\sum_{p \geq J} y'_{i,p} + \sum_{p \geq k} y_{i,p} & \text{if } k \leq J, \\
\sum_{p \geq k} y'_{i,p} & \text{if } k > J,
\end{cases}
\\ \label{eq:product_phi}
\sum_{p \leq k} (y_{i,p} + y'_{i,p}) & =
\begin{cases}
\sum_{p \leq k} y_{i,p} & \text{if } k < J, \\
\sum_{p \leq k} y'_{i,p} + \sum_{p \leq J} y_{i,p} & \text{if } k \geq J.
\end{cases}
\end{align}
\end{subequations}
Additionally, note that
\begin{subequations}
\label{eq:inner_product_sums}
\begin{align*}
\varphi_i(m') - \varepsilon_i(m') & = \inner{h_i}{\wt(m')} = \sum_{p \geq J} y'_{i,p},
\\ \varphi_i(m) - \varepsilon_i(m)  & = \inner{h_i}{\wt(m)} = \sum_{p \leq J} y_{i,p}.
\end{align*}
\end{subequations}

From Equation~\eqref{eq:product_epsilon}, either some $k > J$ obtains the minimum in $\varepsilon_i(m \cdot m')$, in which case $\varepsilon_i(m \cdot m') = \varepsilon_i(m')$, or $\varepsilon_i(m \cdot m') = \varepsilon_i(m) - \inner{h_i}{\wt(m')}$.
Hence, we have
\[
\varepsilon_i(m \otimes m') = \max\bigl(\varepsilon_i(m'), \varepsilon_i(m) - \inner{h_i}{\wt(m')}\bigr) = \varepsilon_i(m \cdot m').
\]
Next, from Equation~\eqref{eq:product_phi}, either some $k < J$ obtains the maximum in $\varphi_i(m \cdot m')$, in which case $\varphi_i(m \cdot m') = \varphi_i(m)$, or $\varphi_i(m \cdot m') = \varphi_i(m') + \inner{h_i}{\wt(m)}$.
Hence, we have
\[
\varphi_i(m \otimes m') = \max\bigl(\varphi_i(m), \varphi_i(m') + \inner{h_i}{\wt(m)}\bigr) = \varphi_i(m \cdot m').
\]
It is clear that $\wt(m \otimes m') = \wt(m) + \wt(m') = \wt(m \cdot m')$.
Thus it remains to show $\psi$ commutes with the crystal operators.

Suppose $k_e(m \cdot m') \leq J$, which is equivalent to $k_e(m \cdot m') = k_e(m)$ since our conditions on $y'_{i,k}$ imply $k_e(m') \neq J$ and
\[
-\sum_{p \geq k} (y_{i,p} + y'_{i,p}) = -\inner{h_i}{\wt(m')} - \sum_{p \geq k} y_{i,p}
\]
for all $k \leq J$.
Thus, we must have
\[
\varphi_i(m') - \varepsilon_i(m') - \varepsilon_i(m) = \sum_{p \geq J} y'_{i,p} + \sum_{p \geq k_e(m)} y_{i,p} < \sum_{p \geq k_e(m')} y'_{i,p} = -\varepsilon_i(m')
\]
from Equation~\eqref{eq:product_epsilon} as $k_e(m \cdot m')$ is the maximum index that achieves the minimal partial sum $-\varepsilon_i(m \cdot m')$. Hence, we have $\varphi_i(m') < \varepsilon_i(m)$, and so
\[
e_i\bigl(\psi(m \otimes m')\bigr) = \psi\bigl(e_i(m \otimes m')\bigr) = \psi(e_i m \otimes m').
\]
Similarly, if $k_e(m \cdot m') > J$, which is equivalent to $k_e(m \cdot m') = k_e(m')$,
then we have
\[
-\varepsilon_i(m \cdot m') = \sum_{p \geq k_e(m')} y'_{i,p} = -\varepsilon_i(m') \leq \sum_{p \geq J} y'_{i,p} + \sum_{p \geq k_e(m)} y_{i,p} = \varphi_i(m') - \varepsilon_i(m') - \varepsilon_i(m)
\]
from Equation~\eqref{eq:product_epsilon} since $k_e(m \cdot m')$ is the maximal index such that $-\varepsilon_i(m \cdot m')$ is achieved.
Thus we have $\varphi_i(m') \geq \varepsilon_i(m)$. Therefore, we have
\[
e_i\bigl(\psi(m \otimes m')\bigr) = \psi\bigl(e_i(m \otimes m')\bigr) = \psi(m \otimes e_i m').
\]

Suppose $k_f(m \cdot m') < J$, which is equivalent to $k_f(m \cdot m') = k_f(m)$ since our conditions on $y_{i,k}$ imply $k_f(m) \neq J$ and
\[
\sum_{p \leq k} (y_{i,p} + y'_{i,p}) = \inner{h_i}{\wt(m)} + \sum_{p \leq k} y'_{i,p}
\]
for all $k \geq J$. Thus we must have
\[
\varphi_i(m') + \varphi_i(m) - \varepsilon_i(m) = \sum_{p \leq k_f(m')} y'_{i,p} + \sum_{p \leq J} y_{i,p} \leq \varphi_i(m \cdot m') = \varphi_i(m)
\]
from Equation~\eqref{eq:product_phi} and that $k_f(m \cdot m')$ is the minimal index that achieves the maximal partial sum $\varphi_i(m \cdot m')$.
Thus, we have $\varphi_i(m') \leq \varepsilon_i(m)$, and so we have
\[
f_i\bigl(\psi(m \otimes m')\bigr) = \psi\bigl(f_i(m \otimes m')\bigr) = \psi(f_i m \otimes m').
\]
Similarly, if $k_f(m \cdot m') \geq J$, which is equivalent to $k_f(m \cdot m') = k_f(m')$, then we have
\[
\varphi_i(m \cdot m') = \varphi_i(m') + \varphi_i(m) - \varepsilon_i(m) > \varphi_i(m)
\]
as $k_f(m \cdot m')$ is the minimal index such that $\varphi_i(m \cdot m')$ is achieved from Equation~\eqref{eq:product_phi}. Hence, we have $\varphi_i(m') > \varepsilon_i(m)$. Therefore, we have
\[
f_i\bigl(\psi(m \otimes m')\bigr) = \psi\bigl(f_i(m \otimes m')\bigr) = \psi(m \otimes f_i m').
\]
\qed\end{proof}

\begin{remark}
Note that the proof of Lemma~\ref{lemma:product_isomorphism} holds for any type (\textit{i.e.}, not just type $A_n^{(1)}$).
Furthermore, Lemma~\ref{lemma:product_isomorphism} could be considered as a generalization of the results of~\cite[Sec.~2]{KS14} and in type $A_n$ is a special case of~\cite[Thm.~2.2]{KTWWY15} (see also~\cite[Prop.~2.6]{KTWWY15}) and the results of~\cite{AN18}.
\end{remark}

\begin{theorem}
\label{thm:nakajima_tensor_product}
Let $j_1, j_2, \dotsc, j_N \in \ZZ_{\geq 0}$ be pairwise distinct. We have
\[
\prod_{k=1}^N \tau_{j_k}(\mcM^{1,s_k}) \iso \bigotimes_{k=1}^N B^{1,s_k}.
\]
\end{theorem}

\begin{proof}
First note that for any permutation $\sigma$ of $\{1, \dotsc, N\}$, we obtain the same crystal
\[
\prod_{k=1}^N \tau_{j_k}(\mcM^{1,s_k}) = \prod_{k=1}^N \tau_{j_{\sigma(k)}}(\mcM^{1,s_{\sigma(k)}})
\]
by the commutativity of variables $Y_{i,k}$ (\textit{i.e.}, the product can be taken in any order).
Next, from the combinatorial $R$-matrix, we have that
\[
\bigotimes_{k=1}^N B^{1,s_k} \iso \bigotimes_{k=1}^N B^{1,s_{\sigma(k)}}
\]
for any permutation $\sigma$ of $\{1, \dotsc, N\}$.
Therefore, we can assume without loss of generality that $j_1 < j_2 < \cdots < j_N$.

Let $\Phi \colon \bigotimes_{k=1}^N B^{1,s_k} \to \prod_{k=1}^N \tau_{j_k}(\mcM^{1,s_k})$ be the map 
\[
\Phi(b_1 \otimes \cdots \otimes b_N) = \prod_{k=1}^N \tau_{j_k}\bigl(\Phi_{s_k}(b_k)\bigr),
\]
where $\Phi_{s_k} \colon B^{1,s_k} \to \mcM^{1,s_k}$ is the isomorphism given by Theorem~\ref{thm:single_row_isomorphism}.
We show the claim holds by induction on $N$. Theorem~\ref{thm:single_row_isomorphism} says this holds when $N = 1$. Thus assume the claim holds for $N-1$.

Consider the tensor product $B \otimes B^{1,s_N}$, where $B = \bigotimes_{k=1}^{N-1} B^{1,s_k}$, and $j_N > j_{N-1}$. From Equation~\eqref{eq:ep_phi_M1s}, we have $-\varepsilon_i(m)$ and $\varphi_i(m)$ being equal to the powers of $Y_{i,1}$ and $Y_{i,0}$, respectively, appearing in $m \in \mcM^{1,s_k}$. Let $\psi \colon \Phi(B) \otimes \tau_{j_N}(\mcM^{1,s_N}) \to \Phi(B) \cdot \tau_{j_N}(\mcM^{1,s_N})$.
Let
\begin{align*}
m & = \prod_{i \in I} \prod_{k \in \ZZ} Y_{i,k}^{y_{i,k}} \in \Phi(B),
\\ m' & = \prod_{i \in I} \prod_{k \in \ZZ} Y_{i,k}^{y'_{i,k}} \in \tau_{j_N}(\mcM^{1,s_N}).
\end{align*}
Note that $\psi(m \otimes m') = m \cdot m'$. Also, for all $i \in I$, we have $y_{i,p} = 0$ for all $p > j_N$, $y'_{i,p} = 0$ for all $p < j_N$ or $p > j_N + 1$ by the definition of $\tau_{j_k}$ and $\mcM^{1,s_k}$. Furthermore, note that $y_{i,j_N}$ is determined by $y_{i,k}$ for all $k < j_N$ for any monomial in $\mcM$.
Hence, $\psi$ is a crystal isomorphism by Lemma~\ref{lemma:product_isomorphism} with $J = j_N$.
\qed\end{proof}

\begin{theorem}
\label{thm:general_tensor_product}
Let $0 < j_1 < j_2 < \dotsm < j_N$ and $\lambda \in P^+$. We have
\[
\prod_{k=1}^{N-1} \tau_{j_k}(\mcM^{1,s_k}) \cdot \tau_{j_N}\bigl( \mcM(\lambda) \bigr) \iso \bigotimes_{k=1}^{N-1} B^{1,s_k} \otimes B(\lambda).
\]
\end{theorem}

\begin{proof}
The proof is similar to the proof of Theorem~\ref{thm:nakajima_tensor_product} with Lemma~\ref{lemma:positivity} implying we can use Lemma~\ref{lemma:product_isomorphism} with $J = j_N$.
\qed\end{proof}

\begin{theorem}
\label{thm:monomial_kyoto_path}
Let $\lambda = \sum_{i \in I} \lambda_i \Lambda_i \in P^+$ be a level $s$ weight, and let $\mu = \sum_{i \in I} \lambda_i \Lambda_{i-1}$.
The map
\[
\Phi \colon \mcM^{1,s} \otimes \mcM(\mu) \to \mcM(\lambda)
\]
given by $m \otimes m' \mapsto m \cdot \tau_1(m')$ is a $U_q'(\asl_n)$-crystal isomorphism.
\end{theorem}

\begin{proof}
First, we note that Lemma~\ref{lemma:positivity} and Lemma~\ref{lemma:product_isomorphism} imply that $m \otimes m' \mapsto m \cdot \tau_1(m')$ is a crystal isomorphism $\mcM^{1,s} \otimes \mcM(\mu) \to \mcM^{1,s} \cdot \tau_1\bigl( \mcM(\mu) \bigr)$.
Therefore, it remains to show that $\mcM^{1,s} \cdot \tau_1\bigl( \mcM(\mu) \bigr) = \mcM(\lambda)$.
Write $\mu = \Lambda_{i_1} + \cdots + \Lambda_{i_s}$, so that $Y_{\mu} = Y_{i_1,0} Y_{i_2,0} \dotsm Y_{i_s,0}$.
Recall from Equation~\eqref{eq:X_variables} that $X_{i,k} = Y_{i-1,k+1}^{-1} Y_{i,k}$.
From~\cite[Cor.~4.9]{Kim05}, we have that a monomial $m \in \mcM(\mu)$ if and only if $m$ can be written as
\[
Y_{j_1,N-1} Y_{j_2,N-1} \cdots Y_{j_s,N-1} \prod_{k=1}^{N-1} X_{1,k-1}^{x_{1,k}} X_{2,k-1}^{x_{2,k}} \cdots X_{n,k-1}^{x_{n,k}},
\]
for some $N \in \ZZ_{\geq 0}$ such that
\begin{enumerate}[($\mu$.1)]
\item $x_{i,k} \in \ZZ_{\geq 0}$ with $\sum_{i=1}^n x_{i,k} = s$ for all $i = 1, 2, \dotsc, n$ and $1 \leq k < N$ and
\item we have $j_p \equiv i_p - (N-1) \pmod{n}$ for all $1 \leq p \leq s$.\footnote{This is called the ground-state condition in~\cite{Kim05}.}
\end{enumerate}
Under $\Phi$, we have
\begin{gather*}
X_{1,0}^{x_{1,0}} X_{2,0}^{x_{2,0}} \cdots X_{n,0}^{x_{1,0}} \otimes Y_{j_1,N-1} Y_{j_2,N-1} \cdots Y_{j_s,N-1} \prod_{k=1}^{N-1} X_{1,k-1}^{x_{1,k}} X_{2,k-1}^{x_{2,k}} \cdots X_{n,k-1}^{x_{n,k}}
\\ \mapsto Y_{j_1,N} Y_{j_2,N} \cdots Y_{j_s,N} \prod_{k=0}^{N-1} X_{1,k}^{x_{1,k}} X_{2,k}^{x_{2,k}} \cdots X_{n,k}^{x_{n,k}}
\end{gather*}
From~\cite[Cor.~4.9]{Kim05}, we need to show that $\Phi(m)$ satisfies
\begin{enumerate}[($\lambda$.1)]
\item $x_{i,k} \in \ZZ_{\geq 0}$ with $\sum_{i=1}^n x_{i,k} = s$ for all $i = 1, 2, \dotsc, n$ and $0 \leq k < N$ and
\item we have $j_p \equiv i'_p + N \equiv i_p + 1 + N \pmod{n}$ for all $1 \leq p \leq s$,
\end{enumerate}
as $Y_{\lambda} = Y_{i'_1,0} Y_{i'_2,0} \dotsm Y_{i'_s,0}$ with $i'_p \equiv i_p + 1 \pmod{n}$ for all $1 \leq p \leq s$ by the definition of $\mu$.
We note that ($\lambda$.1) holds for $k > 0$ since ($\mu$.1) is true and ($\lambda$.1) holds for $k = 0$ from Equation~\eqref{eq:image_X} and the vector realization of $B^{1,s}$.
It is clear that ($\mu$.2) implies ($\lambda$.2).
Hence $\mcM^{1,s} \cdot \tau_1\bigl( \mcM(\mu) \bigr) \subseteq \mcM(\lambda)$, and  $\mcM^{1,s} \cdot \tau_1\bigl( \mcM(\mu) \bigr) \supseteq \mcM(\lambda)$ follows by reversing the above argument.
\qed\end{proof}

\begin{remark}
\label{rem:kyoto_monomial}
Theorem~\ref{thm:monomial_kyoto_path} gives an alternative proof of the existence of the isomorphism given by Equation~\eqref{eq:kyoto_prime} since we have
\[
\Phi_{\lambda}(u_{\lambda}) = \prod_{i \in I} Y_{0,i}^{\lambda_i} = \prod_{i \in I} Y_{0,i}^{\lambda_i} Y_{1,i}^{-\lambda_{i-1}} \times \prod_{i \in I} Y_{1,i}^{\lambda_{i-1}} = \Phi\bigl( \Phi_s(b^{\lambda}) \otimes \Phi_{\mu}(u_{\mu}) \bigr),
\]
where $\Phi_s \colon B^{1,s} \to \mcM^{1,s}$ and $\Phi_{\nu} \colon B(\nu) \to \mcM(\nu)$ are the isomorphisms given by Theorem~\ref{thm:single_row_isomorphism} and Theorem~\ref{thm:highest_weight_monomials} respectively.
\end{remark}

\begin{example}
Consider type $U_q'(\asl_5)$. The ground-state path for $B(\Lambda_0)$ is given by
\[
\young(5) \otimes \young(4) \otimes \young(3) \otimes \young(2) \otimes \young(1) \otimes \young(5) \otimes \young(4) \otimes \young(3) \otimes \young(2) \otimes \young(1) \otimes \cdots.
\]
Therefore, by iterating the isomorphism $\Phi$ from Theorem~\ref{thm:monomial_kyoto_path}, we have
\begin{align*}
\Phi(u_{\Lambda_0}) & = Y_{0,0} Y_{4,1}^{-1} \tau_1(Y_{4,0} Y_{3,1}^{-1}) \tau_2(Y_{3,0} Y_{2,1}^{-1}) \tau_3(Y_{2,0} Y_{1,1}^{-1}) \tau_4(Y_{1,0} Y_{0,1}^{-1}) \tau_5(Y_{0,0} Y_{4,1}^{-1}) \cdots
\\ & = Y_{0,0} Y_{4,1}^{-1} Y_{4,1} Y_{3,2}^{-1} Y_{3,2} Y_{2,3}^{-1} Y_{2,3} Y_{1,4}^{-1} Y_{1,4} Y_{0,5}^{-1} Y_{0,5} Y_{4,6}^{-1} \cdots
\\ & = Y_{0,0}
\end{align*}
If instead we only iterated $\Phi$ five times, we obtain
\[
\Phi(u_{\Lambda_0}) = Y_{0,0} Y_{4,1}^{-1} Y_{4,1} Y_{3,2}^{-1} Y_{3,2} Y_{2,3}^{-1} Y_{2,3} Y_{1,4}^{-1} Y_{1,4} Y_{0,5}^{-1} Y_{0,5} = Y_{0,0}.
\]
\end{example}

By restricting the tensor product given above to a finite number of factors and not including the highest weight crystal, we have
\[
\mcM(Y_{-m,m}^{-s} Y_{0,0}^s) \iso \bigl(B^{1,s})^{\otimes m}.
\]

\begin{remark}
\label{remark:infinity_case}
Our results can be extended to the coherent limit $B_{\infty}$ of $\{ B^{1,s} \}_{s=1}^{\infty}$~\cite{KKM94}. In particular, we can describe $B_{\infty}$ using Nakajima monomials through an extension of Equation~\eqref{eq:image_X}. Then by using the modified Nakajima monomials of Kang, Kim, and Shin~\cite{KKS07} and the analogous characterizations of the monomials appearing in the $B(\infty)$-crystal~\cite[Cor.~4.3]{KKS07}, we can obtain the analogous results of Theorem~\ref{thm:monomial_kyoto_path} for $B(\infty) \to B_{\infty} \otimes B(\infty)$.
\end{remark}

\section{Extensions and Connections}
\label{sec:extensions}

Recall from the introduction that there is a fair amount of evidence that Nakajima monomials should be able to be used to give uniform and natural model for KR crystals. Indeed, $q,t$-characters are given as the sum over Nakajima monomials graded by energy, which admit a classical crystal structure. So the entire $q,t$-character should be constructed by adding $0$-arrows, possibly with the energy grading, from some initial monomial to get a tensor product of KR crystals. In this section, we give some (potential) extensions of our results to more general KR crystals and $q,t$-characters, as well as connections of our results to other crystal models.

\subsection{Single columns}

We can only do $B^{1,s}$ in type $A_n^{(1)}$ using Kashiwara's crystal structure without a quotient precisely because the Dynkin diagram contains a directed cycle. We never repeat an $f_i$ along any shortest path from $b \in B^{1,s}$ back to itself, which matches a path in the Dynkin diagram. So if we express a monomial in terms of the $Y$ variables, we never create $Y_{i,k}$ with $k \notin \{0, 1\}$. This is false in all other orientations of the Dynkin diagram except in the opposite cyclic orientation and for $B^{n,1}$ by duality. However, we are able to construct $B^{r,1}$ using a modification of Nakajima monomials.

We note that elements of $B^{r,1}$ can be described by $(x_1, \dotsc, x_n)$ such that $0 \leq x_i \leq 1$ and $\sum_{i=1}^n x_i = r$ with the crystal structure the same as the vector representation given in Section~\ref{sec:KR_crystals}. Hence, we have the following fact, which does not appear in the literature as far as we are aware, but is likely known to experts.

\begin{proposition}
\label{prop:column_embedding}
There exists a crystal embedding $B^{r,1} \to B^{1,r}$.
\end{proposition}

Therefore, we can quotient our monomials by $X_{i,k}^2$ and obtain a description for $B^{r,1}$ in terms of Nakajima monomials. Explicitly, given a monomial $m$, define the modified crystal operator $\overline{f}_i$ by $\overline{f}_i(m) = f_i(m)$ if $f_i(m)$ does not contain an $X_{i,k}^2$ for some $(i, k) \in I \times \ZZ$ and $\overline{f}_i(m) = 0$ otherwise. The modified crystal operator $\overline{e}_i$ is defined similarly by replacing $f_i$ with $e_i$. Let $\overline{\mcM}(m)$ denote the closure of $m$ under $\overline{e}_i$ and $\overline{f}_i$.

\begin{proposition}
\label{prop:monomial_columns}
For any $k \in \ZZ$, we have
\[
B^{r,1} \iso \overline{\mcM}\left( \prod_{i=1}^r X_{i,k} \right)
\]
as $U_q'(\asl_n)$-crystals.
\end{proposition}

Recall that the highest weight $U_q(\fsl_n)$-crystal $B(\clfw_r)$ corresponding to $B^{r,1}$ can be constructed from the exterior product of the basic representation $B(\clfw_1)$ of $U_q(\fsl_n)$. Therefore, it is more natural to consider the variables $X_{i,k}$ as anticommuting variables to describe $B(\clfw_r)$ using Nakajima monomials (up to a sign).

\subsection{Relations to other models}

Recall the abacus model from~\cite{Tingley08}, which is based on abaci of~\cite[Ch.~2.7]{JK81}. Let $\lambda$ be an level $\ell$ weight. We model $B(\lambda)$ by $\ell$ strings with beads, and the crystal operators act by moving a certain bead.
Fix some abacus configuration $\psi = \{p_j^k \mid j \in I, k \in \ZZ_{\geq0} \}$, where $p_j^k$ denotes the position of the $k$-th bead from the right on the $j$-th string. Let $\psi_i^k = \lvert \{ p_j^k \mid p_j^k \equiv i \mod{n} \} \rvert$. Define
\[
\Gamma(\psi) = \prod_{i \in I} \prod_{k=0}^{\infty} X_{i+1,k}^{\psi_i^k},
\]
where we consider $X_{0,k} = X_{n,k}$.
We note that $\Gamma$ is just a translation of the map $J$ described in~\cite[Thm.~5.1]{Tingley08} in terms of our Nakajima monomial model. Indeed, for each $k$, we can apply Theorem~\ref{thm:single_row_isomorphism}, and by Theorem~\ref{thm:monomial_kyoto_path}, $\Gamma$ is a $U_q'(\asl_n)$-crystal isomorphism.

\begin{example}
Consider the abacus configuration for the highest weight element in $B(2\Lambda_0 + \Lambda_2)$ in type $\asl_5$ 
\[
\begin{tikzpicture}[scale=0.5]
\draw (-7,1) node {$\cdots$};
\draw (-7,0) node {$\cdots$};
\draw (-7,-1) node {$\cdots$};
\foreach \i in {-6,...,0} {
	\fill[black] (\i,1) circle (0.3);
	\fill[black] (\i,0) circle (0.3);
	\fill[black] (\i,-1) circle (0.3);
}
\foreach \i in {1,...,6} {
	\draw (\i,1) circle (0.3);
	\draw (\i,0) circle (0.3);
	\draw (\i,-1) circle (0.3);
}
\fill[black] (1,-1) circle (0.3);
\fill[black] (2,-1) circle (0.3);
\draw (7,1) node {$\cdots$};
\draw (7,0) node {$\cdots$};
\draw (7,-1) node {$\cdots$};
\foreach \i in {-1,0,1}
	\draw[red,-] (5*\i+0.5,1.4) -- (5*\i+0.5,-1.4);
\foreach \i in {-7,...,5}
	\draw (\i+1,-1.8) node {\scriptsize $\i$};
\end{tikzpicture}
\]
where we have drawn the beads as the filled circles and marked the positions of the beads. We have
\[
\psi_{4-k}^k = 2,
\qquad\qquad
\psi_{1-k}^k = 1,
\qquad\qquad
\psi_{-k}^k = \psi_{2-k}^k = \psi_{3-k}^k = 0,
\]
for all $k \in \ZZ_{\geq 0}$. Therefore, we have
\begin{align*}
 \Gamma(u_{2\Lambda_0 + \Lambda_2}) & = \prod_{k=0}^{\infty} X_{2-k,k} X_{-k,k}^2
 \\ & = (Y_{1,1}^{-1} Y_{2,0} Y_{4,1}^{-2} Y_{0,0}^2)
   (Y_{0,2}^{-1} Y_{1,1} Y_{3,2}^{-2} Y_{4,1}^2)
   (Y_{4,3}^{-1} Y_{0,2} Y_{2,3}^{-2} Y_{3,2}^2)
   \cdots
  \\ & = Y_{0,0}^2 Y_{2,0} = \Phi_{2\Lambda_0 + \Lambda_2}(u_{2\Lambda_0 + \Lambda_2}).
\end{align*}
\end{example}

Using the abacus model, we can easily compute the rank-level duality by instead stacking each section of the lift of $\ZZ/n\ZZ$ vertically instead of horizontally (\textit{e.g.}, see~\cite[Fig.~5]{Tingley08}). Therefore, we can compute the rank-level duality on $\mcM(\lambda)$ by conjugating the rank-level duality on the abacus model by $\Gamma$.

Additionally, a bijection $\Xi$ between Nakajima monomials and quiver varieties was constructed in~\cite[Thm.~8.5]{ST14}.
We recall that bijections between the ribbon\footnote{We add this adjective here based on its distinguishing feature: that the crystal operators act by adding a ribbon instead of a single box.} partition model, the abacus model, and cylindric plane partition model for $B(\lambda)$ are given in~\cite{Tingley08}.
Therefore, by using $\Gamma$ and $\Phi$, we can provide another chain of explicit (in the sense of not using the crystal structure) crystal isomorphisms between the quiver varieties and the other three models studied in~\cite{Tingley08}. We remark that such a connection was implicitly already present in the literature. At the end of~\cite[Sec.~2.1]{Tingley08}, a description of how to relate the ribbon partition model with the (multi)partition model via the abacus model, and in~\cite[Thm.~6.4]{ST14}, an isomorphism was given between the (multi)partitions and quiver varieties.

Next, we realize $\widehat{B}^{1,1}$ by the quiver variety represented by an infinite string of arrows
\[
[b] =
\begin{array}{ccccccc}
\cdots & \longrightarrow & (b-2) & \longrightarrow & (b-1) & \longrightarrow & b
\end{array}.
\]
where $f_i [b]$ adds an arrow $b \to b+1$ if $b \equiv i \mod n$ and is $0$ otherwise. We can consider this as an infinite number of bins with a single ball, where the crystal operator $f_i$ moves the ball if it is in the $b$-th bin, where $b \equiv i \mod n$. Thus, we can model $B^{1,1}$ by considering the bins in a cycle of length $n$. Furthermore, from the projection of $\widehat{B}^{1,1}$ onto $B^{1,1}$, we expect a description of $B^{1,1}$ using cyclic quiver varieties~\cite{Nakajima01II} and to have an analogous isomorphism to $\Xi$ with $\mcM^{1,1}$.

\begin{figure}
\begin{center}
\begin{tikzpicture}[scale=0.5]
  \foreach \x in {0, 60, 120, 180, 240, 300}
  {
    \draw (\x:5)+(-0.5,0.5) -- +(-0.5,-0.5) -- +(0.5, -0.5) -- +(0.5, 0.5);
    \draw[->] (\x+10:5) to[in=\x-40,out=\x+100] (\x+50:5);
  }
  \draw (0:6) node{$0$};
  \draw (60:5)+(1,0) node{$1$};
  \draw (120:5)+(-1,0) node{$2$};
  \draw (180:5)+(-1,0) node{$3$};
  \draw (240:5)+(-1,0) node{$4$};
  \draw (300:5)+(1,0) node{$5$};
  \draw[fill=blue!50] (0:5) circle[radius=0.35];
  \draw[fill=red!50] (180:5) circle[radius=0.35];
  \draw[fill=yellow!50] (240:5) circle[radius=0.35];
  \draw[fill=green!50] (266:5) circle[radius=0.35];
  \draw (0:5) node {\small $1$};
  \draw (180:5) node {\small $2$};
  \draw (240:5) node {\small $3$};
  \draw (266:5) node {\small $4$};
  \draw[->, thick, dashed] (250:5.8) to[in=200,out=340] node[midway,below] {$f_3$} (290:5.8);
\end{tikzpicture}
\end{center}
\caption{The model of $B^{1,4}$ for $U_q(\asl_6)$ using balls in bins.}
\label{fig:ball_representation}
\end{figure}

We can extend this ball-bin model into $B^{1,s}$ by using $s$-colored balls, where $f_i$ moves the largest colored ball in bin $i$ (which will be $0$ if there are no balls in bin $i$). For $f_0$, we also need to increase all of the colors of the balls by $1$ and take the result mod $s$. See Figure~\ref{fig:ball_representation} for an example. One can easily check that this gives a crystal isomorphism with the vector realization of $B^{1,s}$ by $x_i$ equals the number of balls in bin $i$.
Similarly, we can construct $\widehat{B}^{1,s}$ by starting with an $s$-fold stack of $[0]$ and $f_i$ adding an arrow to $[b]$ for the minimal $b$ such that $b \equiv i \pmod{n}$ (which will result in $0$ if no such $b$ exists) and $f_i$ adds the arrow on the bottommost such string.

\begin{example}
Consider $\widehat{B}^{1,3}$ for $U_q'(\asl_4)$. We can represent the element $(1, 0, 2, 0)(12)$ by the infinite strings
\[
[13, 15, 15] =
\begin{array}{ccccccccccc}
\cdots & \longrightarrow & 11 & \longrightarrow & 12 & \longrightarrow & 13
\\ \cdots & \longrightarrow & 11 & \longrightarrow & 12 & \longrightarrow & 13 & \longrightarrow & 14 & \longrightarrow & 15
\\ \cdots & \longrightarrow & 11 & \longrightarrow & 12 & \longrightarrow & 13 & \longrightarrow & 14 & \longrightarrow & 15
\end{array}.
\]
Then, we have
\begin{align*}
f_1 [13, 15, 15] & = [14, 15, 15] = 
\begin{array}{ccccccccccc}
\cdots & \longrightarrow & 11 & \longrightarrow & 12 & \longrightarrow & 13 & \longrightarrow & 14
\\ \cdots & \longrightarrow & 11 & \longrightarrow & 12 & \longrightarrow & 13 & \longrightarrow & 14 & \longrightarrow & 15
\\ \cdots & \longrightarrow & 11 & \longrightarrow & 12 & \longrightarrow & 13 & \longrightarrow & 14 & \longrightarrow & 15
\end{array},
\\ f_3 [13, 15, 15] & = [13, 15, 16] =
\begin{array}{ccccccccccccc}
\cdots & \longrightarrow & 11 & \longrightarrow & 12 & \longrightarrow & 13 & \longrightarrow & 14
\\ \cdots & \longrightarrow & 11 & \longrightarrow & 12 & \longrightarrow & 13 & \longrightarrow & 14 & \longrightarrow & 15
\\ \cdots & \longrightarrow & 11 & \longrightarrow & 12 & \longrightarrow & 13 & \longrightarrow & 14 & \longrightarrow & 15 & \longrightarrow 16
\end{array}.
\\
\end{align*}
\end{example}

From Proposition~\ref{prop:monomial_columns}, we can also realize $B^{r,1}$ by considering the $n$ bins placed on a cycle, but now each bin can hold at most one ball.

\subsection{\texorpdfstring{$R$}{R}-matrix kernel}

We describe how our model is a crystal interpretation for the fusion construction of~\cite{KKMMNN91,KKMMNN92}. We note that this could potentially be generalized to give a uniform construction of KR crystals.

From Theorem~\ref{thm:nakajima_tensor_product}, we can construct the tensor product $B^{1,1} \otimes B^{1,1}$ by  considering the product $\mcM^{1,1}$ with its shifted version $\tau_1(\mcM^{1,1})$. However, if we want to avoid the shift, we can still construct $\mcM^{1,1} \otimes \mcM^{1,1}$ by multiplication but having multiplication twisted by a generic parameter $t$. This is a special case of the results of~\cite{Hernandez04,KN12,Nakajima03II,Nakajima03,Nakajima04} expressed in terms of Kashiwara's variation of Nakajima monomials. The kernel of the appropriate parameterized $R$-matrix is generated by
\[
K = \{m \otimes m' - t m' \otimes m \mid m \neq m' \in \mcM^{1,1}\}.
\]
Note that we can construct the elements of $K$ by taking all $t$-commutators\footnote{Recall that a $q$-commutator is $[x,y]_q = xy - qyx$.} for the $q,t$-character, where we consider $\otimes$ as multiplication.

We also recall there is a statistic called \defn{energy} on tensor products of KR crystals~\cite{KKMMNN91,KKMMNN92}. We denote the energy of $b$ by $E(b)$. If we take the quotient of $\mcM^{1,1} \otimes \mcM^{1,1}$ where we consider elements as $t^{E(b)} b$, the quotient by $K$ results in $\mcM^{1,2}$.  This can also be extended to $\mcM^{1,s}$.

\begin{figure}
\[
\begin{tikzpicture}[>=latex,line join=bevel,scale=0.7, every node/.style={scale=0.8}]
\node (node_8) at (236.5bp,301.5bp) [draw,draw=none] {$Y_{0,1}^{-1} Y_{1,0}  \otimes Y_{0,1}^{-1} Y_{1,0} $};
  \node (node_7) at (166.5bp,228.5bp) [draw,draw=none] {$Y_{0,1}^{-1} Y_{1,0}  \otimes Y_{1,1}^{-1} Y_{2,0} $};
  \node (node_6) at (205.5bp,82.5bp) [draw,draw=none] {$Y_{1,1}^{-1} Y_{2,0}  \otimes Y_{0,0} Y_{2,1}^{-1} $};
  \node (node_5) at (356.5bp,301.5bp) [draw,draw=none] {$Y_{1,1}^{-1} Y_{2,0}  \otimes Y_{0,1}^{-1} Y_{1,0} $};
  \node (node_4) at (258.5bp,155.5bp) [draw,draw=none] {$Y_{0,0} Y_{2,1}^{-1}  \otimes Y_{1,1}^{-1} Y_{2,0} $};
  \node (node_3) at (291.5bp,9.5bp) [draw,draw=none] {$Y_{0,0} Y_{2,1}^{-1}  \otimes Y_{0,0} Y_{2,1}^{-1} $};
  \node (node_2) at (288.5bp,228.5bp) [draw,draw=none] {$Y_{0,0} Y_{2,1}^{-1}  \otimes Y_{0,1}^{-1} Y_{1,0} $};
  \node (node_1) at (151.5bp,155.5bp) [draw,draw=none] {$Y_{0,1}^{-1} Y_{1,0}  \otimes Y_{0,0} Y_{2,1}^{-1} $};
  \node (node_0) at (70.5bp,125.5bp) [draw,draw=none] {$Y_{1,1}^{-1} Y_{2,0}  \otimes Y_{1,1}^{-1} Y_{2,0} $};
  \draw [red,->] (node_7) ..controls (162.39bp,208.04bp) and (158.46bp,189.45bp)  .. (node_1);
  \definecolor{strokecol}{rgb}{0.0,0.0,0.0};
  \pgfsetstrokecolor{strokecol}
  \draw (170.0bp,192.0bp) node {$2$};
  \draw [black,<-] (node_2) ..controls (302.46bp,197.19bp) and (308.79bp,180.5bp)  .. (311.5bp,165.0bp) .. controls (321.15bp,109.85bp) and (302.12bp,42.873bp)  .. (node_3);
  \draw (323.0bp,119.0bp) node {$0$};
  \draw [blue,->] (node_7) ..controls (120.66bp,206.65bp) and (82.591bp,194.49bp)  .. (node_0);
  \draw (90.0bp,202.0bp) node {$1$};
  \draw [red,->] (node_5) ..controls (337.15bp,280.3bp) and (317.55bp,259.83bp)  .. (node_2);
  \draw (340.0bp,265.0bp) node {$2$};
  \draw [black,<-] (node_8) ..controls (259.46bp,269.15bp) and (273.86bp,249.49bp)  .. (node_2);
  \draw (278.0bp,265.0bp) node {$0$};
  \draw [red,->] (node_6) ..controls (230.31bp,61.02bp) and (255.96bp,39.844bp)  .. (node_3);
  \draw (268.0bp,46.0bp) node {$2$};
  \draw [red,->] (node_0) ..controls (99.964bp,112.57bp) and (144.0bp,100.62bp)  .. (node_6);
  \draw (140.0bp,109.0bp) node {$2$};
  \draw [blue,->] (node_1) ..controls (166.7bp,134.51bp) and (181.86bp,114.59bp)  .. (node_6);
  \draw (194.0bp,119.0bp) node {$1$};
  \draw [blue,->] (node_8) ..controls (216.58bp,280.3bp) and (196.41bp,259.83bp)  .. (node_7);
  \draw (219.0bp,265.0bp) node {$1$};
  \draw [black,<-] (node_7) ..controls (204.59bp,198.11bp) and (232.19bp,176.8bp)  .. (node_4);
  \draw (232.0bp,192.0bp) node {$0$};
  \draw [black,<-] (node_5) ..controls (368.51bp,243.96bp) and (379.55bp,157.97bp)  .. (335.5bp,110.0bp) .. controls (324.09bp,97.573bp) and (283.16bp,90.69bp)  .. (node_6);
  \draw (377.0bp,192.0bp) node {$0$};
  \draw [blue,->] (node_2) ..controls (280.19bp,207.83bp) and (272.11bp,188.71bp)  .. (node_4);
  \draw (286.0bp,192.0bp) node {$1$};
\end{tikzpicture}
\quad
\begin{tikzpicture}[>=latex,line join=bevel,scale=0.7, every node/.style={scale=0.8}]
\node (node_5) at (99.025bp,9.5bp) [draw,draw=none] {$Y_{0,0}^{2} Y_{2,1}^{-2} $};
  \node (node_4) at (49.025bp,155.5bp) [draw,draw=none] {$Y_{1,1}^{-2} Y_{2,0}^{2} $};
  \node (node_3) at (50.025bp,82.5bp) [draw,draw=none] {$Y_{0,0} Y_{1,1}^{-1} Y_{2,0} Y_{2,1}^{-1} $};
  \node (node_2) at (99.025bp,301.5bp) [draw,draw=none] {$Y_{0,1}^{-2} Y_{1,0}^{2} $};
  \node (node_1) at (127.03bp,155.5bp) [draw,draw=none] {$Y_{0,0} Y_{0,1}^{-1} Y_{1,0} Y_{2,1}^{-1} $};
  \node (node_0) at (50.025bp,228.5bp) [draw,draw=none] {$Y_{0,1}^{-1} Y_{1,0} Y_{1,1}^{-1} Y_{2,0} $};
  \draw [black,<-] (node_0) ..controls (21.915bp,200.07bp) and (8.1749bp,183.1bp)  .. (2.0252bp,165.0bp) .. controls (-7.7459bp,136.24bp) and (21.672bp,106.48bp)  .. (node_3);
  \definecolor{strokecol}{rgb}{0.0,0.0,0.0};
  \pgfsetstrokecolor{strokecol}
  \draw (10.525bp,155.5bp) node {$0$};
  \draw [black,<-] (node_1) ..controls (117.3bp,104.47bp) and (105.18bp,42.128bp)  .. (node_5);
  \draw (123.53bp,82.5bp) node {$0$};
  \draw [blue,->] (node_1) ..controls (109.87bp,140.82bp) and (101.87bp,134.2bp)  .. (95.025bp,128.0bp) .. controls (84.922bp,118.85bp) and (74.046bp,108.12bp)  .. (node_3);
  \draw (103.53bp,119.0bp) node {$1$};
  \draw [black,<-] (node_2) ..controls (108.75bp,250.47bp) and (120.88bp,188.13bp)  .. (node_1);
  \draw (123.53bp,228.5bp) node {$0$};
  \draw [blue,->] (node_2) ..controls (85.302bp,280.62bp) and (71.74bp,260.96bp)  .. (node_0);
  \draw (89.525bp,265.0bp) node {$1$};
  \draw [blue,->] (node_0) ..controls (49.751bp,208.04bp) and (49.489bp,189.45bp)  .. (node_4);
  \draw (58.525bp,192.0bp) node {$1$};
  \draw [red,->] (node_3) ..controls (63.748bp,61.616bp) and (77.311bp,41.964bp)  .. (node_5);
  \draw (89.525bp,46.0bp) node {$2$};
  \draw [red,->] (node_0) ..controls (67.314bp,209.59bp) and (81.83bp,194.95bp)  .. (95.025bp,183.0bp) .. controls (99.196bp,179.22bp) and (103.8bp,175.29bp)  .. (node_1);
  \draw (103.53bp,192.0bp) node {$2$};
  \draw [red,->] (node_4) ..controls (49.299bp,135.04bp) and (49.561bp,116.45bp)  .. (node_3);
  \draw (58.525bp,119.0bp) node {$2$};
\end{tikzpicture}
\]
\caption{The $U_q'(\asl_3)$-crystal $\mcM^{1,1} \otimes \mcM^{1,1}$ (left) and $\mcM^{1,2}$ (right).}
\label{fig:tensor_product}
\end{figure}

\begin{example}
Consider the $U_q'(\asl_3)$-crystals $\mcM^{1,1} \otimes \mcM^{1,1}$ and $\mcM^{1,2}$ (see Figure~\ref{fig:tensor_product}). The graded $q$-character of $\mathcal{T} = \mcM^{1,1} \otimes \mcM^{1,1}$ is
\begin{align*}
\sum_{m \otimes m' \in \mathcal{T}} t^{E(m \otimes m')} m \cdot m' & =
Y_{0,1}^{-2} Y_{1,0}^{2} + \left(t + 1\right) Y_{0,0} Y_{1,1}^{-1} Y_{2,0} Y_{2,1}^{-1} + \left(t + 1\right) Y_{0,0} Y_{0,1}^{-1} Y_{1,0} Y_{2,1}^{-1}
\\ & \hspace{10pt} + \left(t + 1\right) Y_{0,1}^{-1} Y_{1,0} Y_{1,1}^{-1} Y_{2,0} + Y_{1,1}^{-2} Y_{2,0}^{2} + Y_{0,0}^{2} Y_{2,1}^{-2}.
\end{align*}
By considering the (graded) decomposition into $U_q(\fsl_3)$-crystals, we get the same decomposition (after $t \mapsto t^2$) as computed by~\cite{Hernandez04,Nakajima03II,Nakajima03,Nakajima04} (note that we are using different Nakajima monomials). This correspondence is an example of the results of~\cite{KN12}. However, if we take the quotient of the graded $q$-character of the ideal analog of $K$ by sending $t m' \otimes m = 0$ for all $m \neq m' \in \mcM^{1,1}$ (including replacing $\otimes$ with $\cdot$), we obtain
\[
Y_{1,1}^{-2} Y_{2,0}^{2} + Y_{0,0} Y_{1,1}^{-1} Y_{2,0} Y_{2,1}^{-1} + Y_{0,0} Y_{0,1}^{-1} Y_{1,0} Y_{2,1}^{-1} + Y_{0,1}^{-2} Y_{1,0}^{2} + Y_{0,0}^{2} Y_{2,1}^{-2} + Y_{0,1}^{-1} Y_{1,0} Y_{1,1}^{-1} Y_{2,0},
\]
which is the graded $q$-character of $\mcM^{1,2}$.
\end{example}

Our construction works because $B^{1,1}$ in type $A_n^{(1)}$ ``follows'' our orientation of the Dynkin diagram. More generally, in order to construct $\mcM^{1,1}$, we need to take a quotient of a level-zero crystal $\mcM(m)$, where $m$ is some monomial whose weight is of level $0$. For example, in type $C_2^{(1)}$ with $c_{ij} = 1$ if $i < j$ and $0$ otherwise, we can construct $\mcM^{1,1}$ as a quotient of $\mcM(Y_{0,1}^{-1} Y_{1,0})$ by an automorphism $\kappa$ given by $Y_{i,k}^{-1} Y_{i',k'} \mapsto Y_{i,k-2}^{-1} Y_{i',k'-2}^{-1}$ as in Figure~\ref{fig:type_C2}.  However, $\mcM^{1,1} \cdot \mcM^{1,1}$ is not isomorphic to $B^{1,2}$ as the former has 10 elements and the latter has 11. Moreover, $\mcM^{1,1} \cdot \mcM^{1,1}$ agrees with instead taking $\mcM(Y_{0,1}^{-2} Y_{1,0}^2)$ and then taking the quotient by $\kappa$.

\begin{figure}
\begin{center}
\begin{tikzpicture}[>=latex,line join=bevel,xscale=1.1, yscale=0.9, every node/.style={scale=0.9},baseline=0]
  \node (node_8) at (24.0bp,9.5bp) [draw,draw=none] {$Y_{0,4} Y_{1,4}^{-1} $};
  \node (node_5) at (24.0bp,82.5bp) [draw,draw=none] {$Y_{1,3} Y_{2,3}^{-1} $};
  \node (node_4) at (24.0bp,155.5bp) [draw,draw=none] {$Y_{1,3}^{-1} Y_{2,2} $};
  \node (node_3) at (24.0bp,228.5bp) [draw,draw=none] {$Y_{0,3}^{-1} Y_{1,2} $};
  \node (node_6) at (104.0bp,82.5bp) [draw,draw=none] {$Y_{1,1} Y_{2,1}^{-1} $};
  \node (node_2) at (104.0bp,9.5bp) [draw,draw=none] {$Y_{0,2} Y_{1,2}^{-1} $};
  \node (node_1) at (104.0bp,228.5bp) [draw,draw=none] {$Y_{0,1}^{-1} Y_{1,0} $};
  \node (node_0) at (104.0bp,155.5bp) [draw,draw=none] {$Y_{1,1}^{-1} Y_{2,0} $};
  \definecolor{strokecol}{rgb}{0.0,0.0,0.0};
  \pgfsetstrokecolor{strokecol}
  \draw [blue,->] (node_6) ..controls (104.0bp,62.04bp) and (104.0bp,43.45bp)  .. (node_2);
  \draw (109.5bp,46.0bp) node {$1$};
  \draw [red,->] (node_4) ..controls (24.0bp,135.04bp) and (24.0bp,116.45bp)  .. (node_5);
  \draw (32.5bp,119.0bp) node {$2$};
  \draw [black,<-] (node_3) ..controls (78.218bp,186.7bp) and (39.443bp,42.9bp)  .. (node_2);
  \draw (67.5bp,122.0bp) node {$0$};
  \draw [blue,->] (node_5) ..controls (24.0bp,61.722bp) and (24.0bp,42.337bp)  .. (node_8);
  \draw (33.5bp,46.0bp) node {$1$};
  \draw [blue,->] (node_1) ..controls (104.0bp,204.04bp) and (104.0bp,185.45bp)  .. (node_0);
  \draw (112.5bp,192.0bp) node {$1$};
  \draw [blue,->] (node_3) ..controls (24.0bp,213.62bp) and (24.0bp,206.99bp)  .. (node_4);
  \draw (33.5bp,192.0bp) node {$1$};
  \draw [red,->] (node_0) ..controls (104.0bp,131.04bp) and (104.0bp,112.45bp)  .. (node_6);
  \draw (112.5bp,119.0bp) node {$2$};
  \draw [black,->,dashed] (node_8) ..controls (-15.0bp,20.0bp) and (-20.3bp,72.9bp)  .. (-20bp, 100bp);
  \draw [black,<-,dashed] (node_1) ..controls (135.2bp,206.7bp) and (145.4bp,160.9bp)  .. (150bp, 140bp);
%
\end{tikzpicture}
\hspace{60pt}
\begin{tikzpicture}[>=latex,line join=bevel,scale=0.9, every node/.style={scale=1.0}]
  \node (node_2) at (104.0bp,155.5bp) [draw,draw=none] {$Y_{0,2} Y_{1,2}^{-1} $};
  \node (node_6) at (104.0bp,228.5bp) [draw,draw=none] {$Y_{1,1} Y_{2,1}^{-1} $};
  \node (node_1) at (104.0bp,374.5bp) [draw,draw=none] {$Y_{0,1}^{-1} Y_{1,0} $};
  \node (node_0) at (104.0bp,301.5bp) [draw,draw=none] {$Y_{1,1}^{-1} Y_{2,0} $};
  \definecolor{strokecol}{rgb}{0.0,0.0,0.0};
  \pgfsetstrokecolor{strokecol}
  \draw [blue,->] (node_6) ..controls (104.0bp,208.04bp) and (104.0bp,189.45bp)  .. (node_2);
  \draw (112.5bp,192.0bp) node {$1$};
  \draw [red,->] (node_0) ..controls (104.0bp,281.04bp) and (104.0bp,262.45bp)  .. (node_6);
  \draw (112.5bp,265.0bp) node {$2$};
  \draw [blue,->] (node_1) ..controls (104.0bp,354.04bp) and (104.0bp,335.45bp)  .. (node_0);
  \draw (112.5bp,338.0bp) node {$1$};
  \draw [black,<-] (node_1) ..controls (178.218bp,340.0bp) and (178.218bp,190.0bp)  .. (node_2);
  \draw (175.5bp,265.0bp) node {$0$};
\end{tikzpicture}
\end{center}
\caption{A portion of the level-zero crystal $\mcM(Y_{0,1}^{-1} Y_{1,0})$ (left) and $\mcM^{1,1}$ constructed from $\mcM(Y_{0,1}^{-1} Y_{1,0}) / \kappa$ (right) in type $C_2^{(1)}$.}
\label{fig:type_C2}
\end{figure}

Additionally, the quotient of the kernel of the parameterized $R$-matrix corresponding to $B^{1,1} \otimes B^{1,1}$ gives rise to a twisted commutator in type $A_n^{(1)}$. Therefore, we expect that by considering the variables in $\mcM^{1,1}$ as non-commuting variables and then taking an appropriate quotient by the kernel of the $R$-matrix, we could construct general $B^{r,s}$ for all affine types. Furthermore, by relating the extra parameter in the kernel to energy, we expect to recover the results of~\cite{KN12} that relate the $q,t$-characters of standard modules to those of simple modules. Indeed, for the type $C_2^{(1)}$ case considered above, if we add the grading by energy in $\mcM^{1,1} \cdot \mcM^{1,1}$, then we obtain
\begin{align*}
& Y_{0,1}^{-2} Y_{1,0}^2
+ (1 + t) Y_{0,1}^{-1} Y_{1,0} Y_{1,1}^{-1} Y_{2,0}
+ (1 + t) Y_{0,1}^{-1} Y_{1,0} Y_{1,1} Y_{2,1}^{-1}
+ (1 + t) Y_{0,0} Y_{0,1}^{-1}
+ Y_{1,1}^{-2} Y_{2,0}^2
\\ & + (1 + t) Y_{2,0} Y_{2,1}^{-1}
+ (1 + t) Y_{0,0} Y_{1,0}^{-1} Y_{1,1}^{-1} Y_{2,0}
+ Y_{1,1}^2 Y_{2,1}^{-2}
+ (1 + t) Y_{0,0} Y_{1,0}^{-1} Y_{1,1} Y_{2,1}^{-1}
+ Y_{0,0}^2 Y_{1,0}^{-2}.
\end{align*}
However, if we additionally quotient by the kernel of the $R$-matrix, we can construct a Nakajima monomial realization of $B^{1,2}$ as the term $t Y_{0,0} Y_{0,1}^{-1}$ survives but all other terms with a $t$ coefficient are removed.

\subsection{Level-zero construction}

As mentioned in the introduction, the KR crystal $B^{r,1}$ (in general affine type) can be constructed as the quotient of the extremal level-zero crystal $B(\Lambda_r - K^{\vee}_r \Lambda_0)$, where $K^{\vee}_r$ are the dual Kac labels~\cite[Table~Aff1-3]{kac90}, by an automorphism $\eta$ (see, \textit{e.g.},~\cite{NS03,NS05,NS06II,NS06,NS08II}). Here we show how to realize the level-zero crystal in the monomial model (using Kashiwara's version) similar to~\cite{HN06}. However, we do not know how to explicitly determine the automorphism, and hence the quotient, in this model.

\begin{theorem}
\label{thm:extremal_level_zero_monomial}
Let $\varpi = \sum_{i \in I_0} a_i (\Lambda_i - K^{\vee}_i \Lambda_0)$ with $a_i \geq 0$ for all $i$.
Let $\g$ be of affine type and $c_{ij}$ be nonnegative integers for all $i, j \in I$ such that $i \neq j$.
We have
\[
\mcM\left( \prod_{i \in I_0} Y_{0,0}^{-K^{\vee}_i a_i} Y_{i,\ell}^{a_i} \right) \iso B(\varpi)
\]
for any $\ell \geq 0$.
\end{theorem}

\begin{proof}
We refer the reader to~\cite{K95,K02} for the relevant notation (recall that the tensor product convention there is the reverse of ours).
Let $\lambda = \sum_{i \in I_0} a_i \Lambda_i$ and $L = \sum_{i \in I_0} a_i K^{\vee}_i$.
From~\cite[Prop.~8.2]{K95} and~\cite[Eq.~(3.7)]{K02}, we can deduce that the connected component generated by $u_{-L \Lambda_0} \otimes u_{\lambda} \in B(-L \Lambda_0) \otimes B(\lambda)$ is isomorphic to $B(\varpi)$.
%

Let $\mcM$ denote the connected component generated by $Y_{-L \Lambda_0} \otimes \tau_{\ell}(Y_{\lambda}) \in \mcM(-L \Lambda_0) \otimes \mcM(\lambda)$.
From the above and Theorem~\ref{thm:highest_weight_monomials}, we have $\mcM \iso B(\varpi)$. Therefore, it is sufficient to show that there exists a crystal isomorphism $\psi \colon \mcM \to \mcM(Y^{(\ell)}_{\varpi})$, where
\[
Y^{(\ell)}_{\varpi} := \prod_{i \in I_0} Y_{0,0}^{-K^{\vee}_i a_i} Y_{i,\ell}^{a_i}.
\]
Let $\psi$ be the map given by $\psi(m \otimes m') \mapsto m \cdot \tau_{\ell}(m')$.
From Lemma~\ref{lemma:positivity} and its analog for anti-dominant weights, the monomials in $\mcM(\lambda)$ (resp.~$\mcM(-L\Lambda_0)$) only contain $Y_{i,k}$ with $k \geq 0$ (resp.\ $k \leq 0$), and $Y_{i,0}$ must appear with a nonnegative (resp.\ nonpositive) exponent.
Thus, we can apply Lemma~\ref{lemma:product_isomorphism} with $J = 0$ to see that $\psi$ is the restriction of a crystal isomorphism $\mcM(-L\Lambda_0) \otimes \mcM(\lambda) \to \mcM(-L\Lambda_0) \cdot \tau_{\ell}\bigl( \mcM(\lambda) \bigr)$.
Since $\mcM$ is generated by $Y_{-L \Lambda_0} \otimes \tau_{\ell}(Y_{\lambda})$, its image under $\psi$ is generated by $Y^{(\ell)}_{\varpi}$.
Thus $\psi(\mcM) = \mcM(Y^{(\ell)}_{\varpi})$, and hence, $\psi$ is an isomorphism.
\qed\end{proof}

We note that this proof works for any weight $\lambda \in P$ splitting into its dominant part and anti-dominant part (\textit{i.e.}, we write $\lambda^+ + \lambda^- = \lambda$ with $\lambda^+ \in P^+$ and $\lambda^- \in P^-$). Moreover, this construction is likely known to experts.
Furthermore, the proof of Theorem~\ref{thm:extremal_level_zero_monomial} can easily be modified to give a different proof of~\cite[Thm.~2.2]{HN06}.

\appendix
\section{Examples with \textsc{SageMath}}
\label{sec:examples}

We give some examples using {\sc SageMath}~\cite{sage} using the crystal of Nakajima monomials implemented by Ben Salisbury and Arthur Lubovsky.

We construct $B^{1,2}$ for $\asl_5$ using Nakajima monomials and then compare with the tensor product with $B^{1,1}$, verifying Theorem~\ref{thm:nakajima_tensor_product} in this case:
\begin{lstlisting}
sage: P = RootSystem(['A',4,1]).weight_lattice(extended=True)
sage: La = P.fundamental_weights()
sage: c = matrix([[0,1,1,1,0], [0,0,1,1,1], [0,0,0,1,1],
....:             [0,0,0,0,1], [1,0,0,0,0]])
sage: c
[0 1 1 1 0]
[0 0 1 1 1]
[0 0 0 1 1]
[0 0 0 0 1]
[1 0 0 0 0]
sage: M = crystals.NakajimaMonomials(La[1]-La[0], c=c)
sage: x = M({(0,1):-2, (1,0):2}, {})
sage: from sage.categories.loop_crystals import KirillovReshetikhinCrystals
sage: S = M.subcrystal(generators=[x], category=KirillovReshetikhinCrystals())
sage: K = crystals.KirillovReshetikhin(['A',4,1], 1, 2)
sage: K.digraph().is_isomorphic(S.digraph(), edge_labels=True)
True
sage: x = M({(0,1):-2, (1,0):2, (0,2):-1, (1,1):1}, {})
sage: K1 = crystals.KirillovReshetikhin(['A',4,1], 1, 1)
sage: T = tensor([K, K1])
sage: S = M.subcrystal(generators=[x], category=KirillovReshetikhinCrystals())
sage: T.digraph().is_isomorphic(S.digraph(), edge_labels=True)
True
\end{lstlisting}

Next we construct $B^{1,1} \otimes B^{1,2}$ for $\asl_3$:
\begin{lstlisting}
sage: P = RootSystem(['A',2,1]).weight_lattice(extended=True)
sage: La = P.fundamental_weights()
sage: c = matrix([[0,1,0],[0,0,1],[1,0,0]])
sage: M = crystals.NakajimaMonomials(La[1]-La[0], c=c)
sage: x = M({(0,1):-1, (0,2):-2, (1,0):1, (1,1):2}, {})
sage: from sage.categories.loop_crystals import KirillovReshetikhinCrystals
sage: S = M.subcrystal(generators=[x], category=KirillovReshetikhinCrystals())
sage: K1 = crystals.KirillovReshetikhin(['A',2,1], 1,1)
sage: K2 = crystals.KirillovReshetikhin(['A',2,1], 1,2)
sage: T = tensor([K1, K2])
sage: T.digraph().is_isomorphic(S.digraph(), edge_labels=True)
True
\end{lstlisting}

\bibliographystyle{alpha}
\bibliography{monomial}{}
\end{document}